\documentclass[a4paper,11pt]{article}
\usepackage[T2A]{fontenc}
\usepackage[cp1251]{inputenc}
\usepackage[ english]{babel}
\usepackage{amssymb,amsfonts,amsthm,amsmath,mathtext,cite,enumerate,float}
\usepackage[onehalfspacing]{setspace}
\usepackage{geometry}
\usepackage{graphicx}
\usepackage{multirow}
\usepackage{indentfirst}
\usepackage{color}
\usepackage{comment}
\usepackage{datetime}
\date{\displaydate{date}}

\newif\ifbrief
\brieffalse

\ifbrief
\else
\newtheorem{theorem}{Theorem}
\newtheorem{proposition}[theorem]{Proposition}
\newtheorem{corollary}[theorem]{Corollary}
\newtheorem{lemma}[theorem]{Lemma}
\newtheorem{example}{Example}
\newtheorem{remark}[theorem]{Remark}

\theoremstyle{definition}
\newtheorem{definition}{Definition}
\newtheorem{property}[theorem]{Property}

\fi

\author{Ruslan V. Skuratovskii$^{1}$  \\
{\scriptsize  1. Institute of mathematics, NAS of Ukraine,
  ruslan.skuratovskii@imath.kiev.ua, } } 

\title{Normal subgroups of iterated wreath products of symmetric groups}

\date{}
\begin{document}
\maketitle

  \begin{abstract}

Normal subgroups and their properties of finite iterated wreath products of finite symmetric groups and infinite iterated wreath products of symmetric groups are found, the GAP system and theorems of finite group theory are used to identify finite normal subgroups. 
New topology of infinite iterated wreath products of finite symmetric groups are constructed.
All classes of normal subgroups of finite iterated wreath product are investigated, their generators are found and presented in the form of Kaloujnine tables.
Further, the monolith (i.e. the unique minimal normal subgroup) of these wreath products is investigated.
 Conditions when a set-wise stabilizer is a normal subgroup are found.
In the proposed work, the topology of the infinite wreath product as the profinite group is studied, which homeomorphism to the Cantor set is established.
Inverse limit of wreath product of permutation groups is found.


  \end{abstract}

  \section{Introduction}

This research generalizes
our previous results from \cite{SkuNormNPU} and presents a new concept of subwreath product as a subgroup related with normal divisors of the wreath product.

The structure of the normal subgroup lattice of the locally isometry group of the boundary of a spherically homogeneous tree $LIsom\partial T$ was described in \cite{ShLav, Lav, ShLav93}, normality of the stabilizer of the level and rigid stabilizer subgroups were shown in these works. For such kind of infinite groups as $LIsom\partial T$ it was shown that the commutator subgroup belongs to each normal subgroup. The bases of infinite wreath product as one kind of normal divisors were described by Suschansky in \cite{ShNorm}.

   One of the modern branches of group theory research is the study of isometry of geometry structures such as graph trees and of
   combinatorial-topological structures \cite{SkRendi}, in our work the topology of the wreath product as a profinite group is studied, its homeomorphism to the Cantor set is shown.


The main purpose of this study is to describe all normal subgroups of the automorphism group of a finite (restricted) regular rooted
tree.

In \cite{Meld} a structure of commutators of the wreath product $A\wr B$ was briefly considered.
We consider the commutator presentation \cite{Meld} in the form of wreath recursion \cite{Ne, SkRendi } and as Kaloujnine tableaux \cite{LeshCong, LeshPap, Kal, ShNorm}.

Despite this, it was an open question whether the subgroup of $\underset{i=1}{\overset{k}{\mathop{\wr
}}}\,{{{S}}_{{{n}_{i}}}}$ is normal.

  \section{Preliminaries }

Throughout the paper, permutations act on the left.
Let $G$ be a group acting (from the left) by a permutation $\pi$ upon
a set $X$ and the composition of elements will occur from the right to the left, and let $G$  be an arbitrary group. Then the (permutational) wreath product
$G \wr H$ is the semi-direct product $G\ltimes {{H}^{X}}$, where $H$ acts on the direct power ${{H}^{X}}$ by the respective permutations of the direct factors and $|X|=d$. If $a,b \in G \wr H$ which are presented in the form of wreath recursion \cite{Ne}, then its product is given by the formula:
\[ba\,\,=\,\,\sigma ({{b}_{1}},...,{{b}_{d}}) \ \pi\left( {{a}_{1}},...,{{a}_{d}} \right)=\sigma \pi \left( {{b}_{\pi (1)}}{{a}_{(1)}},...,{{b}_{\pi (d)}}{{a}_{(d)}} \right).\]
An inverse element has the form ${{a}^{-1}}={{\pi }^{-1}}\left( a_{{{\pi }^{-1}}(1)}^{-1},...,a_{{{\pi }^{-1}}(d)}^{-1} \right)$.

As is well known, a wreath product of permutation groups is an associative construction.
We will refer to $H^X$ as to the base group of $W$ and to $G$ as the top group of $W$, also $H$ is the bottom group. Let $^{(k-1)}M$ be a direct product $\prod \limits_{i=1}^{k-1}M_{i}$.

According to \cite{ LeshCong, LeshPap, ShNorm, Kal}
 by tableaux we call all possible infinite tuples of the form:
\begin{equation} \label{quo1}
a=[a_k]_{k \in  \mathbb{N}}, \quad\mbox{where}\quad a_k(x_1, x_2, ... , x_{k-1} ) \in
Fun(^{(k)}M,G_{k+1}), \ k\in  \mathbb{N},
\end{equation}
where by $Fun(^{(k)}M, G_{k+1})$ we denote the set of functions from ${}^{(k)}M$ in $G_{k+1}$.

In the current work,
we expand each $a_k$ into a section of the Kalouzhnin table containing $n^{k-1}$ coordinates with elements from $^{(k-1)}M = \prod \limits_{i=1}^{k-1}M_{i}$, provided
we will use the numbering of these sections starting from 0, i.e. $k\in \mathbb{N}\cup\{0\}$.

The $k$-multiple wreath product $\mathop{\wr}\limits^{k}_{i=1}S_{n}$ is denoted by $W_k$ and $\underset{i=1}{\overset{\infty}{\mathop{\wr
}}}\,{{{S}}_{{{n}_{i}}}}$ by $W$.
\begin{definition}\label{k-coordinate}
 We will say that the subset $U < W_k$ is determined by its $k$-coordinate sets $\{U\}_k$, $k \in \mathbb{N}$, if this subset consists of the most possible tables $a \in W_k$ such that for an arbitrary $k \in \mathbb{N}$ the condition $[a]_k \in \{U\}_k$ is satisfied.
\end{definition}

If subsets ${U}_k$ are groups, then a subset $U$ formed by $k$-coordinate sets ${U}_k$ is subgroup of $G$. Each subset $U \subset G$  determined by its $k$-coordinate set ${U}_k$ can be closed to the subgroup of $G$ by taking its group closure which is called a \textit{group formed by $k$-coordinate set} ${U}_k$.

\begin{definition}
    A tableau $u$ has \textbf{depth} $k$ ($d(u)=k$) if all ${{[u]}_{i}}=\varepsilon $ provided $i= \overline{1,2,\,...\,,\,k-1}$; and ${{[u]}_{k}}\ne \varepsilon $. The depth of subgroup $H$ of $W_k$ is $min \{(d(u)),  u \in H \}$ be denoted by $d(H)$.
\end{definition}

Recall that such action of tableau $a=[a_k]_{k \in  \mathbb{N}}$ on an element $m=(m_k)_{k\in \mathbb{N}} \in M$ is defined by the rule on
$m^a=(m_k^{a_k(^{(k-1)}m)})_{k\in   \mathbb{N}}.$

The set $G$ of all such tableaux of the form (\ref{quo1}) with such action which equipped by such multiplication rule forms a
group. The neutral element is the table of the form $e=[\ldots, \varepsilon, \varepsilon, \ldots]$. The trivial subgroup of the finite symmetric group $S_n$ is termed by $E$.

\begin{definition}
 The table $\left[ {{\pi }_{0}},{{\pi }_{1}}\left( {{x}_{0}} \right),...,{{\pi }_{l}}\left( {{x}_{0}},...,{{x}_{l-1}} \right),e,...,e \right]$  of length $l$  be denoted by $^{(l)}u$.
\end{definition}

\begin{definition} We denote by  ${{U}^{(k)}}$ the subgroup of  depth  $\ge k$  and ${}^{(k)}U$ is subgroup of tableaux ${{[u]}_{\alpha }}=\varepsilon$ if $\alpha >k$. An element of $U^{(k)}$ be denoted by $u^{(k)}$ correspondingly.
\end{definition}
In the paper \cite{ShNorm} it was proved that the system of subgroups $U^{(k)}$ forms a normal divisors sequence of $U$.

By $X^n$  we mean vertices at distance $n$ from a root.

 Another convenient language of wreath product presentation is wreath recursion corresponding to the automorphism group of the rooted tree.
 The set $X^*$ is naturally the vertex set of a regular rooted tree, i.e. a connected graph without cycles
and a designated vertex $v_1$ called the root, in which two words are connected by an edge if and only if they are of form $v$
and $vx$, where $v\in X^*$, $x\in X$. Denote by $AutX^*$ the automorphism group of the rooted tree $X^*$.
The set $X^n \subset X^*$ is called the $n$-th level of the tree $X^*$
and $X^1 = \{v_1\}$, another words $X^n :=$ vertices at distance $n$ from root. Denote the $i$-th the vertex of $X^j$ by $v_{ji}$. Note that the unique vertex $v_{ki}$ corresponding to
each word $v$ in alphabet $X$.
For every automorphism $g\in Aut{{X}^{*}}$ and every word $v \in X^{*}$ we define the state $g_{(v)} \in AutX^{*}$ of
$g$ at $v$ by the rule: $g_{(v)}(x) = y$ for $x, y \in X^*$  if and only if $g(vx) = g(v)y = g(v)g_{(v)}(x)$.
The subtree of $X^{*}$ induced by the set of vertices $\cup_{i=1}^k X^i$ is denoted by $X^{[k]}$.

 The restriction of the action of an automorphism $g\in AutX^*$ to the subtree $X^{[l]}$ is denoted by $g_{(v)}|_{X^{[l]}}$. A permutation of automorphism $g\in AutX^*$ acting in a vertex $v_{li}$ is called the \textit{\textbf{vertex permutation}} (v.p.) and denoted by $g_{li}$.

All undeclared terms are from \cite{Ne}.  
The conjugated element $h^{-1}g h$ is denoted by $h^g$.

\begin{definition} A \textit{vertex stabilizer} is the subgroup ${St_{G}{(v)}}=\{g\in G: g(v)\text{ }=v\}$, where $v\in {{X}^{*}}$ is a vertex.
Also $n$-th {level stabilizer} is the subgroup  $S{{t}_{G}}(n)\text{ =}\,\,\bigcap\limits_{v\in {{X}^{n}}}^{{}}{{{G}_{v}}},$ where ${{St}_{G}(v)}$ is the vertex stabilizer.
\end{definition}

We remind a definition of level subgroup $G(l)$ presented by us in \cite{SkCommEur}.

\begin{definition}\label{level subgroup}
We define the subgroup of the $l$-th level $(G(l))$ as a subgroup generated by all possible vertex permutations of this level.
\end{definition}

\begin{definition} \label{coordsubgroup}
  A table of form ${{\left[ a \right]}_{i}}=\varepsilon $ wherein $i\ne k$ be called a $k$-coordinate table.
We denote the subgroup of all $k$-coordinate tables from $W$ by ${{[W]}_{k}}$.
\end{definition}

\begin{remark}\label{Wk}
 The subgroups $W(k)$ and ${{[W]}_{k}}$ are isomorphic presentations of the same subgroup.
\end{remark}
The subgroup ${{[W]}_{k}}$ is also isomorphic to the ${{\mu }_{k-1}}$ -th power of a group ${{S}_{{{n}_{k}}}}$, where ${{\mu }_{k-1}}={{n}_{1}}\times {{n}_{2}}\times ...\times {{n}_{k-1}}$.

Let us fix some notations.
From here and throughout the whole paper we will use the left action in the wreath product.
For convenience we will denote the normal closure of an element $g$ by $N(\widetilde{g})$ if we the group is fixed.

Let $G$ be a group. The intersection of all non-trivial normal subgroups $Mon(G)$ of $G$ is called the monolith of a group $G$.
If  $Mon(G)\neq <1>$, then the group $G$ is called monolithic and, in this case, $Mon(G)$ is the least non-trivial normal
subgroup of $G$.

\begin{definition}
The unique minimal non-trivial normal subgroup is called the monolith.
\end{definition}

Let us introduce some new definitions.
\begin{definition}  A subgroup $U < W = \underset{i=1}{\overset{\infty}{\mathop{\wr
}}}\,{{{S}}_{{{n}_{i}}}}$ is splittable if for each element $u\in U$ and all $l\in \mathbb{N}$ the table ${{u}^{(l)}}$ is also in this subgroup.
\end{definition}
As a corollary we see that in splittable subgroup $U$ the following holds ${}^{\left( l-1 \right)}u \in U$ too.
The splittable normal subgroups of isometry group in the generalized Baire metric space were considered in \cite{Bez}.

In fact, since $u{{=}^{\left( l-1 \right)}}u{{u}^{\left( l \right)}}$ then from existing  ${{u}^{\left( l \right)}}$ for $u\in U$ produces an existing $^{\left( l-1 \right)}u$ in $U$.
This decomposition determines this group structure $U{= ^{(k-1)}}U\rightthreetimes {{U}^{(k)}}$ and additionally ${}^{(k-1)}U \simeq {}^U/{}_{{U}^{(k)}}$.

Now we come to the proof of an important new statement, which was briefly mentioned in the work \cite{Bez} without any proof.

\begin{proposition}\label{k-coordinatesubset}
    If $U$ is a splitting subgroup of the permutational wreath product group $W$, then a group generated by $k$-coordinate tables ${{[U]}_{k}}$ is subgroup of the group $U$, and ${{[U]}_{k}}{{=}^{(k)}}U\bigcap{{{U}^{(k)}}}$.
\end{proposition}
\begin{proof}
 Since splitting of $U$ yields ${{U}^{(k)}} \lhd W$ and ${{U}^{(k+1)}} \lhd {{U}^{(k)}} \lhd U$, we can consider the quotients $ {} [U]_k\simeq {}^{U^{(k)}} /{}_{{U}^{(k+1)}}$ and ${}^{(k)}U \simeq {}^U/{}_{{U}^{(k+1)}}$.

 The remaining part of the proof is elementary by observing that a set is closed under the multiplication operation ${{[U]}_{k}}$ as a quotient group by normal subgroup ${{U}^{(k+1)}}$ and ${{[U]}_{k}}$ is a subgroup in ${{U}^{(k)}}$ because ${{U}^{(k)}}$ splits. This entails ${{[U]}_{k}} < U$.

 In virtue of the isomorphism ${{[U]}_{k}} < {{[W]}_{k}} \simeq \prod \limits_1^{{\mu }_{k-1}} S_{n_k}$, where $ {{\mu }_{k-1}}={{n}_{1}}\times {{n}_{2}}\times ...\times {{n}_{k-1}}$ the following embedding in the group of $k$-coordinate tables is true ${{[U]}_{k}} \hookrightarrow \langle [ \varepsilon,   \ldots, \varepsilon, a_k, \varepsilon, \ldots ]\rangle $. But the group $\langle [ \varepsilon,   \ldots, \varepsilon, a_k, \varepsilon, \ldots ] \rangle$ is a subgroup of $U$.
\end{proof}

\begin{proposition}
A set of open neighborhoods such as the sequence of all $k$-th bases ${{W}^{(k)}}$, where $k\to \infty$ of $\underset{i=1}{\overset{\infty}{\mathop{\wr }}}\,{{{S}}_{{{n}_{i}}}}$, determines a wreath product topology \cite{ShNorm}.
\end{proposition}

  \section{Normal subgroups of wreath product of symmetric groups  }

\subsection{\small Normal subgroups in $S_{n}\wr S_{m}$ }
We will notice at once, that further statements and expositions will be true for $n
\geqslant 3$ even in the case if degrees of two of our symmetric groups from wreath product are different.
 Elements of our wreath product $S_{n} \wr S_{n}, n
\geqslant 5$, will present in the form of tableaux, where $[a]_0$ --- the active element of the table and $a\in S_n$,
$[a_1,a_2,...,a_n]_1$ are the passive elements of the table, $a_i\in S_n$, viz.:
$[a_{}]_{0}, \ [a_{1}, a_{2}, \ldots , a_{n}]_{1}.$
The rule for multiplying elements is described above, as well as in Meldrum's book \cite{Meld}.

Let $k(\pi )$ be the number of cycles in a decomposition of a permutation $\pi $ of degree $n$.
The number $n-k(\pi )$ is  denoted  by $dec(\pi )$, and is called a \emph{decrement} \cite{Sachk} of $\pi $.

As well known \cite{Sachk} the minimal number of transpositions in factorization of a permutation $\pi$ on transpositions is happen to be equal to $dec(\pi)$. We set $dec(e) =0$. Therefore, the decrement of the $n$-cycle is $n-1$.
If $\pi_{1}, \pi_{2} \in S_{n}$, then the following formula holds:
\begin{equation}\label{evendec}
    dec(\pi_{1} \cdot \pi_{2}) = dec(\pi_{1}) + dec(\pi_{2}) -
2m, m \in \mathbb{N},
\end{equation}
where $m$ is number of joint simplifying transpositions in $\pi_1$ and $\pi_2$.
Note, that a decrement of permutation's product $\pi_{1} \pi_{2}$
  can be lesser than a decrement of each permutation $\pi_{1}, \pi_{2}$.

\begin{definition}\label{subwreath product}
	The permutational \emph{ subwreath product }
	$G \wr \hspace{-1mm}\wr H$ is the semidirect product $ G\ltimes {\tilde{H}^{X}}$, where $|X|< \infty$ and $G$ acts on the subdirect product \cite{SubAlg} ${\tilde{H}^{X}}$, that is the kernel of semidirect product, by the respective permutations of the subdirect factors. Provided that the specification of ${\tilde{H}^{X}}$ is established separately.
\end{definition} 

\begin{example}
For instance, in $\mathbb {Z}_n \wr \mathbb{Z}_n$ elements of form:
\[\beta\,\,=\,\,\pi \left( {{b}_{1}},...,{{b}_{n}} \right),\]
where tuples $\left( {{b}_{1}},...,{{b}_{n}} \right)$ satisfying the condition $\mathop{\sum}\limits_{i=1}^{n} b_i \equiv j (mod n)$, $0 \leq j < n$ form a normal subset (set is invariant under conjugation), because of conjugation only permutes elements $a_i$ in the base of $\mathbb {Z}_n \wr \mathbb{Z}_n$; thereby, the sum of residues $\mathop{\sum}\limits_{i=1}^{n} b_i$ remains invariant. Furthermore, the sums of the base elements for the conjugating elements are inverses of each other (by modulo $n$). 

Each value of $j$ corresponds to one normal subset in $\mathbb {Z}_n \wr \mathbb{Z}_n$. Furthermore, if we restrict the set of residues by a set of dividers $d_i, n=p_1^{d_1}p_2^{d_2}...p_k^{d_k} $ of $n$ then we obtain $d(n)$ normal subgroups, where $d(n)=(d_1 +1) \ldots (d_k +1)$ is the number of dividers of $n$ and $p_l$ is a prime factor of $n$, $1\leq l \leq k$. We distinguish normal subgroups of the following type $\mathop{\sum}\limits_{i=1}^{n} b_i \equiv p_l (mod n)$. And the derived classes of normal subgroups of the form $\mathop{\sum}\limits_{i=1}^{n} b_i \equiv 
p_j \times, ... ,\times p_l (mod n),  \, 1 \leq j,l \leq k$ are associated with multiples $p_j, \ldots, p_l$ of these divisors $p_l$.
\end{example}

\begin{example}
A more complex example of a subwreath product of groups that is a normal subgroup $\mathbb {Z}_3 \wr \mathbb{Z}_3$ is a subgroup of elements $a\left( {{b}_{1}},...,{{b}_{n}} \right)$, satisfying the condition $\mathop{\sum}\limits_{i=1}^{n} b_i \equiv a (mod 3)$ and the second $\mathop{\sum}\limits_{i=1}^{n} b_i \equiv 2 a (mod 3)$. These two congruences determine two different normal subgroups that are permutation isomorphic to $C_9 \rtimes C_3$. The element of order 9 in the subgroup satisfying $\mathop{\sum}\limits_{i=1}^{n} b_i \equiv a (mod 3)$ is $[1; 1, 0, 0]$ and for subwreath product determined by $\mathop{\sum}\limits_{i=1}^{n} b_i \equiv 2 a (mod 3)$ an element of order 9 can be generated by $[1; 2, 0, 0]$.


It should be noted that this congruence defines a subwreath product which, in the more general case, appears to be normal subgroups of $\mathbb {Z}_p \wr \mathbb{Z}_p$, $p$ is prime and $p \geq 3$.   
\end{example}

\begin{definition}\label{generalized alternating group}
Let $\widetilde{A}^{(1)}_n$ denote the set of elements of the group $S_{n} \wr S_{n}$ ($n \geqslant 3$) represented as Kaluzhnin tables \cite{Kal} of the form $[e]_{0}, [a_{1}, a_{2}, \ldots , a_{n}]_{1}$ satisfying the condition
\begin{equation}\label{evensum}
\sum_{i=1}^{n}
dec([a_{i}]_{1}) = 2k, \quad k \in \mathbb{N}.
\end{equation}
\end{definition}

In the next proposition, we will prove that $\widetilde{A}^{(1)}_n$ is the normal subgroup of $ S_{n} \wr S_{n}$ so it will be called \textit{\textbf{alternating level
subgroup}} $\widetilde{A}^{(1)}_n$ with  structure $ E \wr\hspace{-1mm}\wr  \widetilde{A}_{n}$, where $\widetilde{A}_{n}$ is the subdirect product uniquely identified by the condition \eqref{evensum}.

It follows directly from the definition that the set of these elements supplemented by the operation of multiplication in the subdirect product, coincides with the group $ E \rtimes (\underbrace{ S_n \boxtimes S_n \boxtimes S_n \boxtimes \ldots \boxtimes S_n}_n )$, where subdirect product satisfies to condition \eqref{evensum}. The subgroup $\widetilde{A}_{n}$ formed by permutations from $X^1$ is called a \textbf{\textit{generalized alternating group}} of permutations from $X^1$. We will identify $\widetilde{A}^{(1)}_n$ and the designation of its structure $E \wr\hspace{-1mm}\wr \widetilde{A}_{n}$.

In other words, such subset of element $g$ from $S_{n} \wr S_{n}, n \geqslant 3$ belongs to class $\widetilde{A}^{(1)}$ if $g$ corresponds to automorphism of $X^{[1]}$ having on zero level the trivial permutation and vertex permutations from $g \in Aut X^{[2]} $ on the first level form an even permutation.

Remind that the improved condition on commutator subgroup of permutation wreath product base was obtained by us in \cite{SkRendi}, in view of the fact that a commutator subgroup is also normal. Now we can conclude that the condition \eqref{evensum} accords with the condition setting the subdirect product in a commutator subgroup of wreath product base \cite{SkRendi}.

\begin{proposition}\label{A2subgr}
$\widetilde{A}^{(1)}_{n}$ forms the normal subgroup in
$ S_{n} \wr S_{n}, n \geqslant 3$.
\end{proposition}

\begin{proof}

  Let us make sure that $\widetilde{A}^{(1)}_n$ is a subgroup.
  Let $\widetilde{a}, \widetilde{b}, \widetilde{c}\in \widetilde{A}^{(1)}_{n} \subset S_n \wr S_n$.
Therefore, an equality of their products $(\widetilde{a} \cdot
\widetilde{b}) \cdot \widetilde{c} = \widetilde{a} \cdot
(\widetilde{b} \cdot \widetilde{c})$ is based on the associativity of $S_n \wr S_n$.
  Since  $\widetilde{A}^{(1)}_n$ is splittable, then by Proposition \ref{k-coordinatesubset} $[\widetilde{A}^{(1)}_n]_1$ is a subgroup of $S_n \wr S_n$.

Besides, it is easy to check this without applying Proposition \ref{k-coordinatesubset},
because for an every element in the form $[e]_{0}, \ [a_{1}, a_{2},
\ldots , a_{n}]_{1}$ there is an inverse element $\widetilde{A}^{(1)}_n$ having the form: $[e]_{0}, \ [a_{1}^{-1}, a_{2}^{-1}, \ldots ,
a_{n}^{-1}]_{1}$.
  Therein, the decrement permanency  $\sum_{i=1}^{n}
dec([a_{i}]_{1}) = \sum_{i=1}^{n} dec([a_{i}^{-1}]_{1}) = 2k$,  $k
\in \mathbb{N}$, $0\leq i \leq n$ holds.

 To prove the normality of $\widetilde{A}^{(1)}_n$ we note that the base $B$ of $S_n \wr S_n$ is isomorphic to $\prod\limits_{i=i}^{n} S_n^{(i)}$ and $\widetilde{A}^{(1)}_n$ has index 2 in $\prod\limits_{i=i}^{n} S_n^{(i)}$.
 Cosets of $B$ by a subgroup $\widetilde{A}^{(1)}_n$ are the following  $\widetilde{A}^{(1)}_n$ and $B\setminus \widetilde{A}^{(1)}_n$ whose elements admit only odd sum of decrements opposite to the condition \ref{evensum}. The last one is not a subgroup because of for two arbitrary elements $\alpha , \beta \in B \setminus \widetilde{A}^{(1)}_n$ the closure of the multiplication operation is violated because of their product $\alpha  \beta \in \widetilde{A}^{(1)}_n$.
 As is well known automorphic image of subgroup is again subgroup and taking into account that $B \lhd S_n \wr S_n$ we obtain that, action of external automorphism $\Theta \in
 S_n \wr S_n$ on $B$ can maps its subgroup $\widetilde{A}^{(1)}_n$ of index 2 only in $\widetilde{A}^{(1)}_n$.
 This confirms the normality of $\widetilde{A}^{(1)}_n$.
\end{proof}

\begin{theorem}\label{A^{(1)}_n}
The order of $\widetilde{A}^{(1)}_{n}$ is $\frac{(n!)^n}{2}$.
\end{theorem}
\begin{proof}

Since $\widetilde{A}^{(1)}_{n} \simeq E \rtimes (\underbrace{ S_n \boxtimes S_n \boxtimes S_n \ldots \boxtimes S_n}_n ) = E \rtimes \widetilde{A}_{n} $ and the second one satisfies the condition of an even sum of decrements \eqref{evensum}. Therefore, permutations with an even sum of decrements in the tuple with index 1 make up exactly half of all permutations contained in the group $E \wr S_n$.

  In view of the fact that $\widetilde{A}_{n}$ admits only a half of the permutations of $\underbrace{ S_n \times S_n \times S_n \ldots \times S_n}_n $, its order is $\frac{(n!)^n}{2}$.
In addition, the subdirect product $\widetilde{A}_{n} \lhd \prod\limits_{i=1}^{n} S_n^{(i)}.$
\end{proof}

\begin{theorem}
The subgroup $\widetilde{A}^{(1)}_3$ of $S_{3} \wr S_{3}$ has the structure
\begin{equation}\label{StrA1}
\widetilde{A}^{(1)}_3 \simeq (C_3 \times C_3 \times C_3) \rtimes (C_2 \times C_2).
\end{equation}
For $n \leq 3$ the structure of the subgroup $\widetilde{A}^{(1)}_n \lhd S_{n} \wr S_{n}$ is $$\widetilde{A}^{(1)}_n \simeq (\prod_{i=1}^{n} A_n ) \rtimes (\prod_{i=1}^{n-1} C_2 ).$$
\end{theorem}
\begin{proof}
To establish the existence of the semidirect product \eqref{StrA1} we consider $Aut(C_3 \times C_3 \times C_3)$.
    Since $Aut(C_3 \times C_3 \times C_3) \simeq GL(3,3)$ has more then 2 involutions due scalar matrices, in addition, its order is $(33-1)\times(33-3)\times(33-32)$,
    therefore homomorphism from $C_2 \times C_2$ in $GL(3,3)$ exists. For further generalization we remark that $C_3 \simeq A_3$ and simple.

For case $n=4$, $Aut(A_4)\simeq S_4$, so the existence homomorphism from $C_2 \times C_2 \times C_2 $ to $Aut(A_4 \times A_4 \times A_4 \times A_4)$ is established.

 In general, for $n>4$ to the structure of $\widetilde{A}^{(1)}_{n} \lhd S_{n} \wr S_{n}, \, n > 3$, the following reasoning is valid.
    Since ${{A}_{n}}$ is an indecomposable and non-abelian group and taking into account that $Aut{{A}_{n}}\simeq {{S}_{n}}$, Theorem 3.1 of J. N. S. Bidwell \cite{Bidw} and moreover, its Corollary 3.3 are applicable. In view of Corollary 3.3 for $n\ge 3$ and $Hom(A_n, Z(A_n))=1$ because of $Z(A_n))=e$, we get $ Aut{{\left( {{A}^n_n} \right)}}=  {{S}_{n}} \wr {{A}_{n}}$.

Since the group $Aut{{\left( {{A}^{n}_{n}} \right)}}\simeq \, {{S}_{n}} \wr {A}_{n}$ contains distinct involutions that can be homomorphic images of generators in $\prod\limits_{i=1}^{n-1}{{{C}_{2}}}$, there exists a homomorphism from $\prod\limits_{i=1}^{n-1}{{{C}_{2}}}$ to $Aut\left( \prod\limits_{i=1}^{n}{{{A}_{n}}} \right)$.
In more detail, as follows from the spectrum of $S_n \wr A_n$, this group contains even more involutions than $n-1$; therefore, it admits a homomorphism from $\prod_{i=1}^{n-1} C_2$.
\end{proof}

Now we can recursively construct the easiest and elegant subgroup $E \wr \widetilde{A}_n^{(2)}$ of $S_{n} \wr S_{n} \wr S_{n}$.

\begin{definition}\label{widetilde{A}_n^{(1)}}
    The subgroup $E  \wr   \widetilde{A}_n^{(1)}$ be denoted by $\widetilde{A}_{n}^{(2)}$.
\end{definition}

Furthermore, we will prove that $ E \wr \widetilde{A}_n^{(1)} \lhd S_{n} \wr S_{n} \wr S_{n}$.

\begin{proposition}
The order of $ \widetilde{A}^{(2)}_n$ is $(n!)^{n^2} \cdot 2^{-n}$.
\end{proposition}
\begin{proof}
    In view of Theorem \ref{A^{(1)}_n} the subgroup $\widetilde{A}^{(1)}_{n}$ has the order $\frac{(n!)^n}{2}$, so we can count the base order of the wreath product $E \wr \widetilde{A}^{(1)}_n$ based on the known order of the bottom group. The construction of a wreath product of permutation groups contains a direct product of a bottom group $\widetilde{A}^{(1)}_{n}$, which determining the order of the base.
\end{proof}

Now we will formulate the equivalent definition of this subgroup.
\begin{definition}\label{A3n}
 The set of elements from $S_{n} \wr S_{n} \wr S_{n}, n
\geqslant 3$ presented by tables \cite{Kal} of the form: $[e]_{0}, \
[e, e, \ldots , e]_{1},  [a_{1}, a_{2}, \ldots , a_{n^2}]_{2},$ satisfying the following condition
\begin{equation}\label{evensum3n}
 \sum_{i=1}^{n}
dec([a_{i}]_{2}) = 2k_1,  \sum_{i=n+1}^{2n}
dec([a_{i}]_{2}) = 2k_2, \ldots,   \sum_{i=n^2- n}^{n^2}
dec([a_{i}]_{2}) = 2k_n, k_j \in \mathbb{N},
\end{equation}
be denoted by
  $\widetilde{A}_{n}^{(2)}$.
\end{definition}

\begin{proposition}
 The definitions \eqref{A3n} and \eqref{widetilde{A}_n^{(1)}} are equivalent.
\end{proposition}
\begin{proof}
    By the definition of base wreath product $E  \wr   \widetilde{A}_n^{(1)}$ contains $n$ examples of group $\widetilde{A}_n^{(1)}$. Each of them  $\widetilde{A}_n^{(1)}$ is the subwreath product described in Definition \ref{evensum} with base $\widetilde{A}_{n} \lhd \prod\limits_{i=1}^{n} S_n^{(i)}$ satisfying the condition \eqref{evensum3n}.
\end{proof}

\begin{proposition}\label{Mon Short}
The subgroup $E \wr E \wr A_n $ is the monolith in $ S_n \wr S_n \wr S_{n}$.
\end{proposition}
\begin{proof}
Now we show that $E \wr E \wr  A_n$ is a minimal and normal.

Consider a conjugation of $\bar{b} \in E \wr E \wr A_n$, $\bar{b}= e ({{\beta}_{1}},...,{{\beta}_{n^2}})$, in ${{S}_{n}}\wr {{S}_{n}}\wr {{S}_{n}} \simeq Aut{{X}^{[3]}}$
\begin{gather}
{{\pi }^{-1}}\left( a_{{{\pi }^{-1}}(1)}^{-1},...,a_{{{\pi }^{-1}}(n^2)}^{-1} \right) e ({\beta}_{1},...,\beta_{n^2})\pi \left( {{a}_{1}},...,{{a}_{n^2}} \right)= \nonumber \label{conjmonolith}
\\[-3mm]
 \label{directconj} 
   \\[-3mm]
\pi^{-1}  {{\pi }} \left( a_{{{\pi }^{-1}} \pi (1)}^{-1}{{\beta}_{\pi (1)}}{{a}_{1}},...,a_{{{\pi }^{-1}} \pi (n^2)}^{-1}{{\beta}_{\pi (n^2)}}{{a}_{n^2}} \right)= e \left( a_{1}^{-1}{{\beta}_{\pi (1)}}{{a}_{1}},...,a_{n^2}^{-1}{{\beta}_{\pi (n^2)}}{{a}_{n^2}} \right),  \nonumber
\end{gather}
    where $\beta_i$ and $a_{ \pi (i)}$ are states of automorphisms $\bar{b} \in {{E}_{n}}\wr {{E}_{n}}\wr {{A}_{n}}$ and $a \in {{S}_{n}}\wr {{S}_{n}}\wr {{S}_{n}}$ in vertices of $X^2$.
    The formula \eqref{directconj} yields that the group $ E \wr E \wr A_n$ is the subgroup of the point-wise stabilizer of $X^2$, this entails the conjugation of its elements by elements of $S_n \wr S_n \wr S_n$ is reduced to coordinate-wise conjugation as in the direct product of $A_n$, perhaps only with subsequent rearrangement of the coordinates themselves under the action by action of $S_n \wr S_n \wr S_n$. Consequently, $E \wr E \wr A_n \lhd S_n \wr S_n \wr S_n$.

    Moreover in view of direct coordinate-wise conjugation established in \eqref{directconj} we have $ E \wr E \wr A_n  \lhd  \widetilde{A}_{n}^{(2)}$ furthermore $ E \wr E \wr A_n  \lhd  \widetilde{A}_{n^2}^{(2)}$ as the subgroup meets the conditions \eqref{evensum3n} and \eqref{widetilde{A}_n^{(1)}}.

    Thus, the following normal series takes place:
    $$E \wr E \wr A_n  \lhd  \widetilde{A}_{n}^{(2)}  \lhd \widetilde{A}_{n^2}^{(2)} \lhd  S_n \wr S_n \wr S_{n}.$$

It should be noted that the smallest normal closure of any nontrivial permutation in each coordinate
from $\mathop{\prod}\limits_{i=1}^{n^2} {A_n}$ in $W_3$ is the direct product of ${{A}_{n}}$.
Since as shown in Galois's prime theorem, the normal closure of any even permutation in $S_n$ is ${{A}_{n}},\,\,n\ne 4$, the normal closure of any even permutation in ${{S}_{n}}$ is ${{A}_{n}}$.

Consider the normal closure of an even permutation ${\pi}_{1}$ in the group $S_n \wr S_n \wr S_n$ on the first coordinate of tuple $e({\beta}_{1},...,\beta_{n^2}) \in E \wr E \wr A_n$ with automorphism states in the vertices of 2-nd level $X^2$
 (the same is for the $k$-th level).

Let $\pi_1$ be 3-cycle. As is well known, the normal closure of a 3-cycle in ${{S}_{n}}$ is ${{A}_{n}}$ since it is generated by all 3-cycles and by conjugation in one 3-cycle one can obtain all 3-cycles.

There is also a second type of even permutations consisting of the product of transpositions $(a, b) (l,k )$ in a first coordinate of $b$, which again reduces to 3 cycles by multiplying two elements $\bar{b}, \bar{b'} \in E \wr E \wr A_n$, containing pairs of transpositions $\left( a,b \right)(l,k )$ and $\left( a,c \right) \left( l, d \right)$ along the first coordinates.

 So if an even permutation is generated by a pair of transpositions, then the problem again reduces to 3 cycles since $\left( a,b \right)\left( a,c \right)=\left( a,b,c \right)$.
 But the product of independent transpositions also reduces to the product of 3 cycles $\left( a,b \right)\left( c,d \right)=\left( a,b \right)\left( a,c \right)\left( a,c \right)\left( c,d \right)=\left( a,b,c \right)\left( a,d,c \right)$ \cite{aviv}. Thus, as is shown above \textbf{normal closure of any non-trivial element} $g \in A_n $ coincides with $A_n$ in view of $A_n$ simplicity.

The similar reasoning can be spread on the rest of coordinates of $(\beta_{1},...,\beta_{n^2})$.
The same reasoning applies for each state $\beta_i \in ({{\beta}_{1}},...,{{\beta}_{n^2}}), \, 1 \leq i \leq n^2$ coordinate from the ${{X}^{2}}$ of automorphism $b$  since the conjugation is coordinate-wise, as in the direct product consequently $e ({{\beta}_{1}},...,{{\beta}_{n^2}}) ^{W_3} \simeq \mathop{\prod}\limits_{i=1}^{n^2} {A_n}$.

Therefore a normal closure of arbitrary non-trivial an element $b= e ({{b}_{1}},...,{{b}_{n^2}})$, where $b_i \in A_n $ are states in vertices of $X^2$, in $W_2$ is subgroup of $E\wr E\wr{A_n}$.

Assume that some split-extension $N$ of $E \wr E \wr A_n$ by any subgroup of $S_n \wr S_n \wr S_n$ is a normal subgroup too, then it has form $ N = E \wr H \wr A_n$ and consequently contains $E \wr E \wr A_n$, thus such $N$ is not minimal.

Any subgroup of the form $E \wr H \wr E$ or $H \wr E \wr E$ has a trivial intersection with the base of $S_n \wr S_n \wr S_n$, therefore, according to the well-known theorem on the non-trivial intersection of any normal subgroup with the base of the wreath product, $E \wr H \wr E$ or $H \wr E \wr E$ are not normal.

Denote the subgroup $E \wr E \wr A_n$ by $M_3$.
Assume that there is another normal subgroup $K$ such as $K \cap{M_3} \neq e$ then $K$ contains an element having an even permutation $\pi$ on $X^2$,
then $\pi \in M$ too and $\pi^{W_3} = \mathop{\prod}\limits_{i=1}^{n^2} {A_n}$ due to conjugation in $W_3$ spread this permutation on all coordinates of the last level $X^2$.

Hence, $M_3 \simeq \mathop{\prod}\limits_{i=1}^{n^2} {A_n}$ is the monolith of $W_3$.
The proof is completed.
\end{proof}
\begin{corollary}\label{Mon Short E wr A_n}
The subgroup $ E \wr A_n $ is the monolith in $ S_n \wr S_{n}$ provided $n>4$ or $n=3$, for the case $n=4$ the monolith is $ E \wr K_4$.
\end{corollary}
\begin{proof}
    Restricting the group structure considered in Proposition \ref{Mon Short} to the group $E \wr A_n$ for both the cases $n>4$, $n=3$ and $n=4$, we obtain by the same reasoning a proof that $E \wr A_n$ is a unique minimal normal subgroup in $S_n \wr S_n$.
\end{proof}

\begin{definition}\label{A^k}
{\sl The set of elements from $\mathop{\wr}\limits^{k}_{i=0}  S_{n_i}, n_i
\geq 3$ with depth $k$
satisfying the following condition}
\begin{gather} \label{evensumk}
 \sum_{i=1}^{n^{k}}
dec([a_{i}]_{k}) = 2t, t \in \mathbb{N}.    \, \
\end{gather}
 be denoted by $\widetilde{A}^{(k)}_{n^{k}}$ and called \textbf{generalized alternating $k$-th level subgroup}. By $\widetilde{A}^{(0)}_{n^{ }}$ we mean the usual $A_{n}$.
\end{definition}
The fact that $\widetilde{A}^{(k)}_{n^{k}}$ having a structure $E \wr \ldots \wr E \, \wr \hspace{-2,5mm}\wr \, \widetilde{A}_{n_k^k}$, where $\widetilde{A}_{n_k^k}$ is a subdirect product uniquely identified by the condition \eqref{evensumk}, and has property of normality will be proved in Corollary \ref{normalAk}. To prove that $\widetilde{A}_{n^k}^{(k)}$ is the normal subgroup we can argue as in Proposition \ref{A2subgr}. Checking of a subgroup conditions is similar to the arguments from the first point of the Proposition \ref{A2subgr} proof.

Now we can recursively construct the easiest and elegant subgroup $E \wr \widetilde{A}_n^{(1)}$ of $S_{n} \wr S_{n} \wr S_{n}$.
\begin{definition}\label{widetilde{A}_n^{(1)}}
    The subgroup with the structure $E  \wr   \widetilde{A}_n^{(1)}$ be denoted by $\widetilde{A}_{n}^{(2)}$.
\end{definition}

 Corollary \ref{normalAk} shows us that $\widetilde{A}_{n}^{(2)} \lhd W_2$. Furthermore, in Remark \ref{series} we show that $\widetilde{A}_{n}^{(2)} \lhd \widetilde{A}_{n^2}^{(2)} \lhd S_{n} \wr S_{n} \wr S_{n}$.

\begin{theorem}
    The maximal normal subgroups of $W_k$ are $M_l = \mathop{\wr}\limits^{l-1}_{i=1}  S_{n_i} \wr \hspace{-2,5 mm} \wr \, \widetilde{A}^{}_{n_l^{l}} \wr (\mathop{\wr}\limits^{k}_{i=l+1} S_{n_i})$, when $1 < l \leq k$ and $M_1 =  \, A_{n_1} \wr (\mathop{\wr}\limits^{k}_{i=2} S_{n_i})$ if $l=1$.
\end{theorem}
\begin{proof}
Firstly we note that in the case $l=1$ the subdirect product $\widetilde{A}^{}_{n^{l}}$ degenerates in $\prod\limits_{i=1}^{n^l} {A}^{}_{n^{}}$.
The base observation in this proof is that the index  $\mid W_k : M_l \mid =2$ for $1 \leq l \leq k$, thence we obtain both maximality and normality for each $1 \leq l \leq k$.
\end{proof}

We generalize Definition \ref{A^k}  on some levels with even sum of permutation decrements.

\begin{definition}\label{Am,k}
 The set of elements from $\mathop{\wr}\limits^{k}_{i=1}     S_{n_i}, n_i
\geqslant 3$ with depth $m$

satisfying the following condition
\begin{equation}\label{evensummtok}
\sum_{i=1}^{n^j}
dec([a_{i}]_{j}) = 2t, t \in \mathbb{N}, \,\,  m \leq j \leq k,   \, \, \
[a_{i}]_{j} = e,  \, whenever \,
j=\overline{ 0,  m-1}
\end{equation}
 be called $\widetilde{A}^{(m, k)}_{n^{j}}$, where $m<k$.
\end{definition}
 Note, if we assume that $m=k$, then this group degenerates into the group $\widetilde{A}^{(k)}_{n^{k}}$.

It is easy to verify axioms of group for $\widetilde{A}^{(m, k)}_{n^{j}}$,
 its normality is shown in Corollary  \ref{{A}^{(k)}_{n^{m}}<W_k}.

\begin{lemma}\label{Invunderconj}

    Let $\alpha, \beta \in W_k$ and $\alpha^{}= \pi^{} ( a_{{{ }^{}}(1)}^{},...,a_{{{ }^{}}(n^j)} )$,  $\alpha^{-1}= \pi^{-1} ( a_{{{\pi}^{-1}}(1)}^{-1},...,a_{{{\pi }^{-1}}(n^j)}^{-1} )$ $\beta^{}= \sigma ( b_{{{ }^{}}(1)}^{},...,b_{{{ }^{}}(n^j)})$, wherein $a_{i}^{}$, $b_{i}^{}$ are states of automorphisms $\alpha$ and $\beta$ respectively in $v_{ji}$, on $j$-th level, $m \leq j \leq k$, then:

    1)  the sums of $a_{{{ }^{}}i}$ decrements are invariant under inversion of elements $\alpha, \alpha^{-1}$:
    $$\sum\limits_{i=1}^{n^j}{dec(a_{{{\pi }^{-1} }(i)}^{-1})}=\sum\limits_{i=1}^{n^j}{dec({{a}_{i}})},$$

    2) for a pair of conjugate elements $\alpha$, $\alpha^{\beta}$ the equality $\sum\limits_{i=1}^{n^j}{dec(ab_{{{\sigma }^{}} }a_{{{\pi } \sigma \pi^{-1}}(i)}^{-1})}=\sum\limits_{i=1}^{n^j}{dec({{b}_{i}})}$ holds.
    \end{lemma}

\begin{proof}
The states ${a}_{i}$ of automorphism $\alpha$ are permuted only under the action of $\pi^{-1}$ in $\alpha^{-1}$, thus the sum of the decrements $\sum\limits_{i=1}^{n^j}{dec(a_{{{\pi }^{-1} }(i)}^{-1})}$ remains constant with respect to the initial $\sum\limits_{i=1}^{n^j}{dec({{a}_{i}})}$.
The special case for $\pi = e$ of the statement from item 1) has already been proven in Proposition \ref{A2subgr}. Thus, the statement of item 1) holds.

Consider a conjugation of $\beta \in S_n \wr S_n$ by an arbitrary element $\alpha$.
\begin{gather}
{{\pi }^{-1}}\left( a_{{{\pi }^{-1}}(1)}^{-1},...,a_{{{\pi }^{-1}}(n^j)}^{-1} \right) \sigma ({{b}_{1}},...,{{b}_{n^j}})\pi \left( {{a}_{1}},...,{{a}_{n^j}} \right)= \nonumber
\\[-3mm]
\label{conjXk}
   \\[-3mm]
\pi^{-1}  \sigma {{\pi }} \left( a_{{{\pi }^{-1}}\sigma \pi (1)}^{-1}{{b}_{\pi (1)}}{{a}_{(1)}},...,a_{{{\pi }^{-1}}\sigma \pi (n^j)}^{-1}{{b}_{\pi (n^j)}}{{a}_{(n^j)}} \right).  \nonumber
\end{gather}

Taking into account the invariance of the decrements of $a \in Aut({X}^{[k]}) $
with respect to taking the inverse element and the fact that the acting of $\pi $ and $\sigma$ only rearranges states $a_i$ and $b^{-1}_i$ of automorphisms $\alpha$ and $\beta$ in the right part of \eqref{conjXk} over the vertices of $X^{j}$, $1 \leq j \leq k$ in the tuple of elements $( b_{{{\pi }^{-1}}\sigma \pi (1)}^{-1}{{a}_{\pi (1)}} {{b}_{(1)}},...,b_{{{\pi }^{-1}}\sigma \pi (n^j)}^{-1}{{a}_{\pi (n^j)}}{{a}_{(n^j)}} )$ respectively to the initial element
$\beta$ and $\alpha$
we conclude that $\sum\limits_{i=1}^{n^j}{dec(ba_{{{\sigma }^{}} }b_{{{\pi^{-1} } \sigma \pi}(i)}^{-1})}=\sum\limits_{i=1}^{n^j}{dec({{a}_{i}})}$, $m \leq j \leq k$, which completes the proof.
\end{proof}
\begin{corollary}\label{normalAk}
The set of elements $\widetilde{A}^{(k)}_{n^{k}}$ is the normal subgroup of  $W_k$, where $k \geq 2$.
\end{corollary}
\begin{proof}
In view of Lemma \ref{Invunderconj} the conjugation by elements of $W_3$ described in \eqref{directconj}, where permutation $\sigma =e$, preserves the level subgroup $\widetilde{A}_{n^k}^{(k)}$ subject to the condition \eqref{evensumk} of Definition \ref{A^k}, thence the requirements of parity \eqref{evensumk} from Definition \ref{A^k} are satisfied after conjugation in $W_k$.
    The closure under multiplication is easily verified in a completely analogous manner as for $\widetilde{A}^{(1)}_{n^{}}$ in Proposition \ref{A2subgr}. The Definitions \ref{A^k} and \ref{subwreath product} entails the structure of this group  $E \wr \ldots \wr E \, \wr \hspace{-2,5mm}\wr \, \widetilde{A}_{n_k^k}$, where $\widetilde{A}_{n_k^k}$ is a subdirect product mentioned in Definition \ref{A^k}.

    Also this statement can be proved as in Proposition \ref{A2subgr} with using the fact that $\widetilde{A}_{n^k}^{(k)}$ has index 2 in the base of $W_k$.
\end{proof}

The study of invariant subgroups that arise as an embedding $\widetilde{A}_{{{n^{}}^{m}}}^{(k)}$ in $\widetilde{A}_{{{n^{}}^{k}}}^{(k)}$, where $m<k \in \mathbb{N}$, will be considered in the following theorem.
\begin{remark} \label{series}
{ Generalized alternating subgroup ${\tilde{A}_{{{n^{}}^{k}}}^{(k)}}$ is not solvable for $n>4$. The following normal series is in place among the subgroups of the $k$-th level:
\begin{equation}\label{serieschain}
E \lhd  \widetilde{A}^{(k)}_{n^{0}}  \lhd \widetilde{A}_{{{n^{}}^{}}}^{(k)} \lhd \widetilde{A}_{{{n^{}}^{2}}}^{(k)} \lhd \widetilde{A}_{{{n^{}}^{3}}}^{(k)} \lhd \ldots
  \widetilde{A}_{{{n^{}}^{k}}}^{(k)} \lhd W_k.
\end{equation}
}
    \end{remark}
\begin{proof}
   The key step of this proof is the observation that elements of $\widetilde{A}^{(k)}_{n^{m}}$ satisfy condition \ref{evensumk} of Definition \ref{A^k} thence $\widetilde{A}^{(k)}_{n^{m}}< \widetilde{A}^{(k)}_{n^{k}}$, for $m<k$.
  In addition, checking the satisfaction of the condition, very similar to \eqref{directconj} for the case $k>2$, which is easy to do due to Lemma \ref{Invunderconj}, indicates the normality of $\widetilde{A}^{(k)}_{n^{m}}$. Furthermore when $m=1$ we have $\widetilde{A}^{(k)}_{n^{}}$ that is the monolith of $W_k$, in addition we clarify $\widetilde{A}^{(k)}_{n^{}} = \prod\limits_{m=1}^k {A}_{n}$. Hence, this subgroup finalize each normal series of generalized alternating groups.

The quotient group for this series for $1<j<k$ the a following subdirect product satisfying the even number of 1 condition,

$${ }^{\tilde{A}_{{{n^{}}^{j+1}}}^{(k)}} / {}_{\tilde{A}_{{{n^{}}^{j}}}^{(k)}} \simeq \mathop{\prod}\limits_{i=1}^{{{n}^{}}} \mathbb{Z}_2. $$

But the last quotient space $\widetilde{A}^{(k)}_{n^{0}}$, appearing for $j=1$, is a monolith of $W_k$, isomorphic to $\prod\limits_{m=1}^k {A}_{n}$, which is not abelian for $n>4$, so ${\tilde{A}_{{{n^{}}^{k}}}^{(k)}}$ is not solvable.
\end{proof}

Recall that in the case of the left action, a diagonal subgroup $Diag(S_n \wr S_m)$ of $S_n \wr S_m$ is the subgroup consisting of functions $g_2(x_1)=const$ having the same values on $X_1$ so $Diag(S_n \wr S_m) \simeq S_m$. In addition, we remark that its elements are represented by Kalouzhnine tables of the form $[1, g_2(x_1)]$.
Obviously $Diag(S_n \wr S_m) \lhd S_n \wr S_m$.
Now we generalize the concept of a diagonal subgroup.

\begin{definition}\label{T_n^{(1)}}
{\sl A subgroup in  $S_{n} \wr
S_{n}$ is called \textbf{generalized diagonal subgroup of $1$-st level } and is denoted by $\widetilde{T}^{(1)}_{n}$ if it consists of:
\begin{enumerate}
\item elements of $ E \wr A_{n}\,$, note that elements of the first item form the subgroup $e \wr A_{n} $,
\item elements with the tableau \cite{Kal} presentation  $[e]_{1}, \ [\pi_{1}, \ldots ,\pi_{n}]_{2}$, where $\pi_{i} \in S_n \setminus A_n$,
\end{enumerate}}
The specification of homomorphism $\phi :{{C}_{2}}\to Aut({{A}_{n}}\times {{A}_{n}}\times ...\times {{A}_{n}})$ is defined by the action of non-trivial element $h\in {{C}_{2}}$ corresponds to automorphism ${{\phi }_{h}}\in Aut({{A}_{n}}\times {{A}_{n}}\times ...{{A}_{n}})$ which consists in left action $h\cdot n=\left( (12)\cdot {{\pi }_{1}},...,(12)\cdot {{\pi }_{n}} \right)$, where ${{\pi }_{i}}$ is permutation on $i$-th coordinate of $n\in {{A}_{n}}\times {{A}_{n}}\times ...\times {{A}_{n}}$.
    Also we note that the operation $\boxplus$ of a subdirect product is subject of item 1) and 2).        

Its structure $E  \wr \hspace{-2,5mm}\wr \, \widetilde{T}^{}_{n}$, where $\widetilde{T}^{}_{n}$ is a subdirect product uniquely identified by the conditions of items 1) and 2).

\end{definition}
The structure of the subdirect product $\widetilde{T}^{}_{n}$ emerging in the definition is as follows:
 $$\tilde{T_n} \simeq  (\mathop{ \underbrace{A_n \times A_n \times \cdots
    \times A_n}  }\limits_{n} )  \rtimes  C_2 \simeq \mathop{ \underbrace{S_n \boxplus S_n \ldots    \boxplus S_n}\limits_{n} },$$
    where the operation $\boxplus$ of a subdirect product is the subject of items 1) and 2).

One can easily validate the correctness of this definition, namely, that the set of such elements forms a subgroup.

\begin{proposition}
    The order of $\widetilde{T}^{(1)}_{n}$ is $\frac{(n!)^n}{2^{n-1}}$.
\end{proposition}
\begin{proof} A tuple of the first type elements from Definition \eqref{T_n^{(1)}} has the size
$\frac{(n!)^n}{2^{n}}$. Further, a tuple of the second type of $\widetilde{T}^{(1)}_{n}$ has the same size. Thence, $\frac{(n!)^n}{2^{n}} + \frac{(n!)^{n}}{2^{n}} = \frac{(n!)^{n}}{2^{n-1}}$.
\end{proof}

\begin{proposition} \label{T^{(1)}_{n}}
    The group $\widetilde{T}^{(1)}_{n}$ is normal in $S_n \wr S_n$.
\end{proposition}
\begin{proof}
Consider a conjugation of $b \in  E \wr \hspace{-2,5mm}\wr \, \widetilde{T}_n$, $b= e ({{b}_{1}},...,{{b}_{n}})$, in $ {{S}_{n}}\wr {{S}_{n}}\wr {{S}_{n}} \simeq Aut{{X}^{[3]}}$
\begin{gather}
{{\pi }^{-1}}\left( a_{{{\pi }^{-1}}(1)}^{-1},...,a_{{{\pi }^{-1}}(n)}^{-1} \right) e ({{b}_{1}},...,{{b}_{n}})\pi \left( {{a}_{1}},...,{{a}_{n}} \right)= \nonumber \label{conjmonolith}
\\[-3mm]
\label{conjT_n}
   \\[-3mm]
\pi^{-1}  {{\pi }} \left( a_{{{\pi }^{-1}} \pi (1)}^{-1}{{b}_{\pi (1)}}{{a}_{(1)}},...,a_{{{\pi }^{-1}} \pi (n)}^{-1}{{b}_{\pi (n)}}{{a}_{(n)}} \right)= e \left( a_{(1)}^{-1}{{b}_{\pi (1)}}{{a}_{(1)}},...,a_{(1)}^{-1}{{b}_{\pi (n)}}{{a}_{(n)}} \right).  \nonumber
\end{gather}
In view of Lemma \eqref{Invunderconj} the equality $dec (a_{ (i)}^{-1}{{b}_{ (i)}}{{a}_{(i)}}) \equiv dec (a_{{{\pi }^{-1}} \pi (i)}^{-1}{{b}_{\pi (i)}}{{a}_{(i)}})$.

This group $\widetilde{T}^{(1)}_{n}$ does not admit any split extension for the following reason. Assume that the complementary subgroup $A$ is non-trivial and permutation $\sigma \in A$.
In order to satisfy Definition \ref{T_n^{(1)}} decrements values of the functions of a first tuple of tableaux have to be the same parity under conjugation by arbitrary $a$ i.e.
{\small
\begin{gather}\label{parity} dec (a_{{{\pi }^{-1}}\sigma \pi (1)}^{-1}{{b}_{\pi (1)}}{{a}_{(1)}}) \equiv dec (a_{{{\pi }^{-1}}\sigma \pi (2)}^{-1}{{b}_{\pi (1)}}{{a}_{(2)}}) \equiv \cdots \equiv dec(a_{{{\pi }^{-1}}\sigma \pi (n^k)}^{-1}{{b}_{\pi (n^k)}}{{a}_{(n^k)}}) mod 2, \end{gather}
}
That implies ${{\pi }^{-1}}\sigma \pi (i) = i$, which is only possible iff $\sigma=e$.

Given that  Definition \ref{T_n^{(1)}} ensures that $$dec ({{b}_{(1)}}) \equiv dec ({{b}_{(2)}}) \equiv \cdots \equiv dec({b}_{(n)}) mod 2, $$
and therefore $dec ({{b}_{\pi(1)}}) \equiv dec ({{b}_{\pi(2)}}) \equiv \cdots \equiv dec({{b}_\pi{ (n)}}) mod 2 $,
we conclude that for condition \ref{parity} to hold, the permutation $\sigma$ must be trivial.
In order for a functions $a_{i}$ of a first tuple of a tableau to be invariant under conjugation they must be invariant under action of $\sigma$.
\end{proof}

We remark that $\widetilde{T}^{(1)}_{n} \ntriangleleft \widetilde{A}^{(1)}_{n}$ for $n\cong 1 \mod n$ because of $\sum\limits_{i=1}^{n^j}
dec([a_{i}]_{1})\cong 1mod2$ hence does not satisfies condition \ref{evensum} of Definition \ref{generalized alternating group} that is indicated on Fig. 1.

\begin{definition}\label{T_n^2}
{\sl A subgroup in  $S_{n} \wr
S_{n} \wr
S_{n}$ is denoted by $\widetilde{T}_{n^2}^{(2)}$ if it consists of: 
\begin{enumerate}
\item elements of the form $E \wr E \wr A_{n}\,$,
\item elements with the tableau \cite{Kal} presentation   $[e]_{0}, \ [e, \ldots ,e]_{1}, \ [\pi_{1}, \ldots ,\pi_{n}, \pi_{n+1}, \ldots ,\pi_{n^2}]_{2}$, wherein $\forall i
= 1, \ldots , n$: $ \pi_{i} \in S_n \setminus A_n$.
\end{enumerate}  }
\end{definition}

We recursively define the subgroup $\widetilde{T}_n^{(2)}$ having $n$ different intervals of elements with the same parity permutations on $X^2$.
\begin{definition}\label{T_n^{(3)}}
The subgroup of $S_{n} \wr S_{n} \wr
S_{n}$ having structure $ E \wr \widetilde{T_{n}} $ be denoted by $\widetilde{T}_{n}^{(2)}$.
\end{definition}

\begin{property}
The following isomorphism $ \widetilde{T}_{n}^{(2)} \simeq \mathop{ \underbrace{S_n \boxplus S_n \ldots    \boxplus S_n}\limits_{n} } \times \mathop{ \underbrace{S_n \boxplus S_n \ldots  \boxplus S_n}\limits_{n} } $ $ \times \ldots   \times  \mathop{ \underbrace {S_n \boxplus S_n \ldots    \boxplus S_n}\limits_{n} }$, where a tuple $S_n \boxplus S_n \ldots    \boxplus S_n$ repeats $n$ times, holds. The operation $\boxplus$ is determined in Definition \ref{T_n^{(1)}}.       
\end{property}

\begin{proof}
The group structure of the wreath product $ E \wr \widetilde{T_{n}}$ provides a $n$-multiple direct product of $\widetilde{T_{n}}$ in its base.
 The operation $\boxplus$ accords with the properties described in the items 1 and 2 of Definition \eqref{T_n^{(1)}}; also, $\boxplus$
 is determined by the automorphism in $\tilde{T_n} \simeq (\mathop{ \underbrace{A_n \times A_n \times \cdots
    \times A_n}  }\limits_{n})  \rtimes  C_2$ in this case.   
\end{proof}

\begin{remark}
Note that in $\widetilde{T}_n^{(2)}$ vertex permutation of a tableau second part $[e]_{0}, \ [e \ldots ,e]_{1},$  $ [\pi_{1} \ldots ,\pi_{n}; \pi_{n+1}
    \ldots ,\pi_{n^2}]_{2}$ satisfies the condition
  that either all $ \pi_{i_j} \in S_n \setminus A_n$ or all $ [\pi_{i_j}]_3 \in A_n$ for $1 \leq  i_1 \leq n$, $n+1 \leq i_2 <2n$, ... , $n^2-n+1 \leq  i_{n} \leq n^2$.
\end{remark}

 \begin{remark} \label{T3}
The subgroups $\widetilde{T}_{n^2}^{(2)}$, $\widetilde{T_{n}}^{(2)}$ of  $W={{S}_{n}}\wr {{S}_{n}}\wr {S}_{n}$ are normal.    
\end{remark}
\begin{proof}
 Since $\widetilde{T}_{{{n}^{2}}}^{(2)}<S{{t}_{W}}\left( {{v}_{0}} \right)$ then via conjugating by automorphism $\pi \left( {{g}_{1}},...,{{g}_{n}} \right)$ in $W$ the automorphisms $\left( {{g}_{1}},...,{{g}_{n}} \right)$ multiplies as in direct product. Furthermore, each vertex permutation  
which is situated on ${{X}^{2}}$ lay in one of $S{{t}_{W}}\left( {{v}_{1i}} \right)$ where $1\le i\le n$ thence the tuples of permutations  ${{g}_{2{{\pi }^{-1}}(1)}},...,{{g}_{2{{\pi }^{-1}}(n)}}$ and ${{g}_{21}},...,{{g}_{2n}}$ in $\pi_{1i}^{-1}\left( {{g}_{2{{\pi }^{-1}}(1)}},...,{{g}_{2{{\pi }^{-1}}(n)}} \right) \times$ $\varepsilon \left( {{h}_{21}},...,{{h}_{2n}} \right){{\pi }_{1i}}\left( {{g}_{21}},...,{{g}_{2n}} \right)$, $\pi_{1i}^{-1}\left( {{g}_{2{{\pi }^{-1}}(n+1)}},...,{{g}_{2{{\pi }^{-1}}(2n)}} \right)\times$ $\varepsilon \left( {{h}_{2,n+1}},...,{{h}_{2,2n}} \right){{\pi }_{1i}}$ $\left( {{g}_{2(n+1)}},...,{{g}_{2,2n}} \right)$ and till the last block 
$\pi _{1i}^{-1}\left( {{g}_{2,{{\pi }^{-1}}({{n}^{2}}-n+1)}},\,...\,,{{g}_{2,\,{{\pi }^{-1}}({{n}^{2}})}} \right)$ \\ $ \varepsilon \left( {{h}_{2,\,{{n}^{2}}-n+1}},\,...\,,{{h}_{2,{{n}^{2}}-n}} \right){{\pi }_{1i}}\left( {{g}_{2,{{n}^{2}}-n+1}},\,...\,,{{g}_{2,{{n}^{2}}}} \right)$ multiplies directly without action of root permutation from $\widetilde{T}_{{{n}^{2}}}^{(2)}$, because it contains stabilizers $St\left( {{v}_{1i}} \right)$, $1\le i\le n$.  The  proof for  ${{\widetilde{T}}_{n}}^{(2)}$ is very similar. 
\end{proof}

\begin{definition}\label{T_{n_k}^{(k)}}
{\sl A subgroup in $\mathop{\wr}\limits^{k}_{i=0}  S_{n}$ is denoted by $\widetilde{T}_{n^k}^{(k)}$ if it is isomorphic to quotient group consisting of two cosets of $\prod\limits_{i=0}^{{{n}^{k}}}{S_{n}^{(i)}}$ by $\prod\limits_{i=0}^{{{n}^{k}}}{A_{n}^{(i)}}$:
\begin{enumerate}
\item elements of $\mathop{\wr}\limits^{k}_{i=0}  E \wr A_{n}\,$,
\item elements from another coset containing only odd permutations in the base of  $\mathop{\wr}\limits^{k}_{i=0}  S_{n}$.
\end{enumerate}  
One can easily validates the correctness of this definition, i.e. that the set of such elements form a subgroup and its normality. This subgroup has the
structure $$\tilde{T_n} \simeq  (\mathop{ \underbrace{A_n \times A_n \times \cdots
    \times A_n}  }\limits_{n^k} )  \rtimes  C_2 \simeq \mathop{ \underbrace{S_n \boxplus S_n \ldots    \boxplus S_n}\limits_{n^k} },$$ }
    \end{definition}

\begin{proposition}
The following homomorphism  $Wr_n / (\mathop{\wr}\limits^{k-1}_{i=1} E \wr A_{n_k}) \cong Wr_{n-1} \wr C_2$ holds.    
\end{proposition}
\begin{proof}
 Construct the homomorphism $Wr_n  \rightarrow Wr_{n-1}\wr C_2$ provided by the element-wise mapping
 $\varphi ([\pi;\rho_1,\rho_2,\ldots,\rho_{p_n}])=   [\pi;\delta_1,\delta_2,\ldots,\delta_{p_n}]$, wherein $\pi \in \mathop{\wr}\limits^{k-1}_{i=1} S_{n_i}$, which is the element-wise mapping defined by the formula ${p_n}=\prod_{i=1}^{k-1} n_i$
  and \\
\begin{equation}
\delta_{n_i}=  \begin{cases}
1,& \text{ if } \rho_{n_i} \in S_{n_k}
\backslash A_{n_k}, \\
0,& \text{if }  \rho_{n_i} \in A_{n_k}.
 \end{cases} \end{equation}
\end{proof}

\begin{definition}\label{A^mk}
The set of elements from $\mathop{\wr}\limits^{k}_{i=1}  S_{n_i}, n_i
\geqslant 3$ with depth $k$
satisfying the following condition
\begin{equation}\label{evensummk}
 \sum_{i=(s-1)n^m+ 1 }^{sn^m}
dec([a_{i}]_{k}) = 2t, t \in \mathbb{N},   \, 1 \leq s \leq n^{k-m}, \, \   
[a_{i}]_{j} = e,  \, for \,
j=\overline{0, k-1 }
\end{equation}
 denoted by $\widetilde{A}^{(k)}_{n^{m}}$, where $m \leq k$.
\end{definition}

Below we show that $\widetilde{A}^{(k)}_{n^{m}}$ forms a normal subgroup with a structure $E \wr \ldots \wr E \wr \hspace{-2mm}\wr \widetilde{A}_{n_k^m}$, where $m<k$. Its special case $\widetilde{A}^{(k)}_{n^{0}} = \prod\limits_{i=1}^{n^k} A_n$ corresponds to the direct power of $A_n$.

\begin{proposition}\label{stab}
   A stabilizer of block $St_{H}(B^s_{n^m})$ fix each block in orbit of a block $B^s_{n^m}$ under action of $H$ iff $St_{H}(B^s_{n^m}) \lhd H$, where $H \lhd W_k$.
\end{proposition}
    \begin{proof}
    The set of $n^m$ vertices of $X^k$ be denoted by $B^s_{n^m}$, where $s$ is correspondent index from Definition \ref{A^mk}.
By $St_{H}(B^s_{n^m})$ we mean set-wise stabilizer of block $B^s_{n^m}$.
    For convenience we call such a block $B^s_{n^m}, \, 1\leq s \leq n^m$ (of vertexes) by a point.

In general, the stabilizers $St_{H}(B^s_{n^m})$ of points from one orbit $O_H (B^1_{n^m})$ are conjugate $St_{H}(B^s_{n^m}) \simeq \alpha St_{H}(B^s_{n^m}) \alpha^{-1}$.
But according to the statement conditions $St_{H}(B^s_{n^m})$ fixes each block from orbit $O_H (B^1_{n^m})$,
 then this stabilizer is common to all points, so we get the equality $St_{H}(B^s_{n^m}) = \alpha St_{H}(B^s_{n^m}) \alpha^{-1}$, which immediately implies that $St_{H}(B^s_{n^m}) \lhd H$.

Vice versa if we have normality i.e. $St_{H}(B^s_{n^m})= \alpha St_{H}(B^s_{n^m}) \alpha^{-1}$, $\alpha \in H$.

Since conjugate stabilizers are stabilizers of different points in an orbit and these stabilizers coincide, this means that the original stabilizer $St_{H}(B^s_{n^m})$ of the block $B^s_{n^m}$ as a set-wise stabilizer fixes all points in the orbit of this block. By points in orbit in this case we mean blocks which are in this orbit as in hyper space where block is a point.
\end{proof}

The previous Proposition \ref{stab} immediately entails normality of a subgroup of $St_{H}(B^j_{n^m})$ having index 2. Stabilizer described in Proposition \ref{stab} we denote by $Rst_{H}(B^j_{n^m})$.

\begin{corollary} \label{{A}^{(k)}_{n^{m}}<W_k}
    The subgroup $\widetilde{A}^{(k)}_{n^{m}}$ is normal in $W_k$.
\end{corollary}
\begin{proof}

Let $H \lhd G$ then we show that if an automorphism $\alpha \in G$ acts on a normal subgroup $N$ of $H$, then its automorphic image $N^\alpha$ is also a normal subgroup of $H$.  Let $g\in G$ then the condition of normality takes form $gNg^{-1} = N$ and $\alpha(g^{-1})=\alpha^{-1}(g)$. Consider its automorphic image $\alpha(g)\alpha(N)\alpha(g^{-1})= \alpha(g)\alpha(N)\alpha^{-1}(g)=\alpha(N)$ that confirms the normality of $N^{\alpha}$.

Since $\widetilde{A}^{(k)}_{n^m}$ is a subgroup of index 2 in a block stabilizer $St_{W_k}(B^s_{n^m})$ and $St_{W_k}(B^s_{n^m}) \lhd W_k$ i.e. $St_{W_k}(B^s_{n^m})$ meets the condition of Proposition \ref{stab}, thence
$\widetilde{A}^{(k)}_{n^m}$ is normal in $St_{H}(B^s_{n^m})$.
Provided this normal subgroup $\widetilde{A}^{(k)}_{n^m}$ is unique with such size in $St_{W_k}(B^s_{n^m})$ then each automorphism in $W_k$ keeps it invariant in particular $g \widetilde{A}^{(k)}_{n^m} g^{-1} = \widetilde{A}^{(k)}_{n^m}$, $g\in W_k$.

 Besides, this corollary can be proved by direct application of Lemma \ref{Invunderconj} to checking the condition \ref{evensummk} of Definition \ref{A^mk} after the conjugation by arbitrary $\alpha \in W_k$.
 \end{proof}

\begin{remark}
  The set $\widetilde{A}^{(k)}_{n^{m}}$ forms subgroup of $\mathop{\wr}\limits^{k}_{i=1}  S_{n_i}$ and having a structure
 $E \wr \ldots \wr E \wr \hspace{-2mm}\wr \widetilde{A}_{n_k^m}$, where $\widetilde{A}_{n_k^m}$ is the subdirect product uniquely identified by the condition \eqref{evensummk}.
\end{remark}
The proof is a direct verification of subgroup conditions, similar to the first part of Proposition \ref{A2subgr} proof.

Consider an infinite wreath product of permutation groups and define a topology on this infinite group. We set $W = \mathop{\wr}\limits^{\infty}_{i=1} S_{n_i}, \,\,  n_i
\geqslant 3$.

Let $\phi_i$ be natural projection $ W_{i+1} \xrightarrow[\phi_i ]{}W_{i}, i\in \mathbb{N} $, defined due to the quotient ${}^{W_{i+1}}/_{W_{i+1}^{(i+1)}}$.
 Then projective limit is $W_{ }=\underset{i \in \mathbb{N}} {\varprojlim}( W_i, \phi_i)= \mathop{\wr}\limits^{\infty}_{i=1} S_{n_i}.$

 \begin{proposition}
    The maximal normal subgroups of $W$ are $$ \widetilde{W}^{k} = S_{n} \wr S_{n} \wr  \ldots \wr  S_{n_{}}  \wr \hspace{-2,5mm} \wr \, \widetilde{A}_{n^k} \wr S_{n_{}} \wr S_{n_{}} \ldots, $$ for $k \geq 2$ and for $k=1$ the maximal normal subgroup of this type is $\widetilde{W}^{1} = A_{n} \wr ( \mathop{\wr}\limits^{l-1}_{i=1}  S_{n}).$

     In the exceptional case $l=1$ the maximal normal subgroup takes the form $\widetilde{W}^{1} = A_{n} \wr ( \mathop{\wr}\limits^{l-1}_{i=1}  S_{n}).$

\end{proposition}
\begin{proof}
      The key observation is $| \widetilde{W}: \widetilde{W}^{k}| = 2$. The minimality of the subgroup index $\widetilde{W}^{k}$ indicates both its normality and its maximality.
   The exceptional case $l=1$ produces the maximal normal subgroup degenerates in $\widetilde{W}^{1} = A_{n} \wr ( \mathop{\wr}\limits^{l-1}_{i=1}  S_{n}).$
   In the case $l=1$ the subdirect product $\widetilde{A}^{}_{n^{l}}$ degenerates in $\widetilde{W}^{1} = A_{n} \wr ( \mathop{\wr}\limits^{l-1}_{i=1}  S_{n})$, it implies that $\widetilde{W}^{1}$ is of index 2.
   \end{proof}

In our paper \cite{SkCommEur}, we introduced the concept of the ${{W}_{k}}(k-1)$ subgroup of $Aut{{X}^{[k]}}$ that consists of v.p. which are located on ${{X}^{k-1}}$. Now we note that ${{\widetilde{W}^{k}}_{}}(k)$ is isomorphic to $ \widetilde{A} ^{(k)}_{n^k}$.

A rigid level stabilizer \cite{Ne} of $l$-th level is termed by $rist_l( \widetilde{W}^{k})$.
\begin{remark}
	The topological group $W_{}$ is homeomorphic to \textbf{Cantor set} $\mathbb{C}_{\frac{1}{3}}$.
\end{remark}

 \begin{proof}
We assign as an open set $O_{n^k}^{(k)}$ the normal subgroup defined by the formula:  $$ O_{n^k_k}^{(k)} = S_{n} \wr S_{n} \wr  S_{n}\wr \ldots \wr  S_{n_{}}  \wr \hspace{-3mm} \wr \, \widetilde{A}_{n^k} \wr S_{n_{}} \wr S_{n_{}} \wr \ldots,$$ where  $1 \leq k \leq \infty$, these subgroups be denoted by $\widetilde{W}^{k}$. The index $[\widetilde{W}^{k} : rist_l( \widetilde{W}^{k})]< \infty$ this yields that $\widetilde{W}^{k}$ is a branch group, moreover it is a just-infinite branch group \cite{GrBranch, GrInfBranch, DueselldorfBranch}.

 In terms of the $k$-level subgroups of $W_{}$ we see that $W_{}(k) \simeq \widetilde{A}^{(k)}_{n^k}$ has a structure of a subdirect product. In a special case $l \neq k$, this subgroup decomposes in the direct product $W_{\infty}(l) \simeq  \prod\limits_1^{n^l} {S}_{n}$.
In terms of $k$-coordinate subsets, we see that $O_{n^k}^{(k)}$ is determined by its $k$-coordinate subsets $\{O_{n^k}^{(k)} \}_k = \widetilde{A}_{n^k}$.

 The subgroup $ S_{n} \wr S_{n} \wr  S_{n}\wr \ldots \wr  S_{n_{}}  \wr \hspace{-3mm} \wr \, \widetilde{A}_{n^m} \wr S_{n_{}} \wr S_{n_{}} \wr \ldots$ of $O_{n^k}^{(k)}$ be denoted by $O_{n^m}^{(k)} $, where $\widetilde{A}_{n^m}$ is from  Definition \eqref{A^mk}, therefore $O_{n^m}^{(k)} \subseteq O_{n^k}^{(k)}$.

 The sets $O_{n^m}^{(k)}$, $1\leq m \leq k$ are declared as open sets of $W_{\infty}$. The boundary case provided $m=1$ is $O_{n}^{(k)} \simeq \prod\limits_{i=1}^k A_n$, wherein each group is both open and closed, that's each point of $O_{n}^{(k)}(k)$ is clopen set. Thus, the topology on $k$-th level subgroup of $\widetilde{W}^{k}(k)$ is discrete for each $k$, as a corollary the separation axioms $T_1$, $T_3$ are satisfied.

To ensure continuity under multiplication, we declare that the cosets of the form $(ij)O_{n^k_k}^{(k)}$, $(ij)\widetilde{O}^{( k)}_{n^{m}}$ are also included in the base of the topology.

 In  $W_{} = \mathop{\wr}\limits^{\infty}_{i=1} S_{n}, \,\,  n_i \geqslant 3$,
 the open set of level $k$ will be denoted as $O_{n^k}^{(k)}$.

If an element $\alpha \in W_{}$ meets Definition \eqref{A^mk} and its $k$-coordinate table $[\alpha]_k \in \widetilde{A}^{( k)}_{n^{m}}$ but $\alpha \notin \widetilde{A}^{(k)}_{n^{m-1}}$ then
a minimal open neighborhood containing an element $\alpha $ is $\widetilde{O}^{( k)}_{n^{m}}$, $m<k$.

For each pair of elements $\alpha ,\beta : \alpha \cdot \beta \in O_{n^k}^{(k)}$ there are $O_{n^m}^{(k)}$ and $O_{n^l}^{(k)}$, $l, m \leq k$ such that $O_{n^m}^{(k)}\ni \alpha $, $O_{n^l}^{(k)} \ni \beta$ and $O_{n^m}^{(k)} O_{n^l}^{(k)} \unlhd O_{n^{max\{l,m\}}}^{(k)}$, provided that it is important that the result of a product of these elements is also in this open neighborhood $\alpha \cdot \beta \in O_{n^{max\{l,m\}}}^{(k)}$.
The coordinate subgroup $[O_{n^m}^{(k)}]_k$ described in Definition \eqref{coordsubgroup} is isomorphic to $\widetilde{A}^{( k)}_{n^{m}}$.

The operation of finding the inverse element $\alpha^{-1}$ for each $\alpha \in \widetilde{O}^{( k)}_{n^{m}}$, $1 \leq m \leq k$ is closed with respect to the same open set $\widetilde{O}^{( k)}_{n^{m}}$, because of by virtue of point 1) of Lemma \ref{Invunderconj} we see that $\alpha^{-1} \in \widetilde{O}^{( k)}_{n^{m}}$. Thus, the quotient ${}^{W_{}} / {}_{O_{n^k_k}^{(k)}}$ has exactly two quotient classes.

The open sets defined on $W_{}(k)$ forms topological space which be denoted by $W^k(k)$.
Tyhonov product of finite sets  $T = \prod\limits_{k=0}^{\infty }{W^k(k)}$ produce cylindrical topology.
Thus, $W_{} \cong \mathbb{C}_{\frac{1}{3}}$ as the set that is a countable product of sets $W^k(k)$ having discrete topology.

According to the well-known corollary 10.60 from Brouwer's theorem \cite{Kainth}, a countable product of finite metric spaces $M(k)$ spaces endowed with discrete topology such that $|M(k)| > 1$ for all $k \in \mathbb{N}$ is homeomorphic to the standard Cantor middle-thirds set
 \cite{Shaw, Kainth, tails}. Thus, $W_{}$ is homeomorphic to every Cantor space $C$ by Moore-Kline Theorem \cite{Shaw, AlgCantor}.
At this point the proof can be fully completed, but we will additionally show a second way to the proof using Brouwer's Theorem.

To show the absence of isolated points, we note that every point $x\in W_{\infty}$ is contained in a certain open neighborhood $\widetilde{O}^{( k)}_{n^{m}}$ which contains a countable direct product of finite sets which, by definition of $\widetilde{O}^{( k)}_{n^{m}}$ contains other points from $W_{\infty}$ so $x$ is not isolated.

Tyhonov product of finite sets  $\prod\limits_{k=0}^{\infty }{{{W}^{k}}(k)}$ produce cylindrical topology. Let ${{x}_{\alpha }}\in {{W}^{\alpha }}(\alpha )$.
In the product topology  (cylindrical topology) each point $\left\{ {{x}_{\alpha }} \right\}_{0}^{\infty }$ is not open set, because open sets are exactly cylindrical set  of this topology. Cylinder sets in products of discrete sets have finite number of fixed coordinates. But an isolated point have to be precisely an open set. Thus this set has no isolated points.
We have a cylindrical topology on Tyhonov product of finite sets ${{W}^{k}}(k)$.

The infinite product space $\prod\limits_{k=0}^{\infty }{{{W}^{k}}(k)}$ has no isolated points because the requirement that $|{{{W}^{k}}(k)}|>1$ for infinitely many indices $k$ prevents any basic open set from collapsing to a single point.
Thus, this space is perfect.

The proof that the Tikhonov product of discrete spaces ${{{W}^{k}}(k)}$ is totally disconnected relies on the fact that the connected component of any point in such a product consists of itself. Consider two points $x=\{{{x}_{\alpha }}\}$ and $y=\{{{y}_{\alpha }}\}$ in $T$. Since $x\ne y$, there are two coordinates such that ${{x}_{i}}\ne {{y}_{i}}$. Including these points in the set of fixed points of the cylindrical set, we obtain two different open-closed neighborhoods that separate $x$ and $y$.
\end{proof}

For the metrization, we introduce the following metric for the space of automorphisms $\alpha ,\beta \in Aut X^{\infty}$ similar to the metric for the space of one-sided sequences
$$\rho \left( \alpha ,\beta  \right)=\sum\limits_{i=0}^{\infty }{\frac{1}{{{2}^{i+1}}}{{\rho }_{i}}\left( {{a}_{ij}},{{b}_{ij}} \right),\,\,\,\,\,1\le j\le {{n}^{i}}},$$
where $a_{ij}, b_{ij}$ are v.p. of $\alpha ,\beta$ from $i$-th level.

Returning to the initially introduced open sets $O_{{{n}^{k}}}^{(k)}$, we analyze the quotient space by the subgroups they support.
In order to show that a finite topology is defined on each factor by $O_{n^k}^{(k)}$, we will show that the number of factor classes over an open set is finite. If the open set is $O_{{{n}^{k}}}^{(k)}$, then there are only 2 factor classes, i.e. 0 and 1. For quotient by $O_{{{n}^{m}}}^{(k)}$ we have $2^{n^{k-m}}$ open sets.

The quotient $W_{\infty} $ by $\widetilde{W}^{k}$ is isomorphic to $\mathbb{Z}_2$ and has the form $ E \times  \ldots \times E  \times \mathbb{Z}_2 \times E \times E \ldots $.
Consider $\cap_{k=1}^{\infty} \widetilde{W}^{k}$ as the kernel of homomorphism from $W_{\infty}$, in terms of normal subgroups this kernel $\cap_{k=1}^{\infty} \widetilde{W}^{k} \simeq \widetilde{A}_{n^k_k}^{(0,\infty)}$.
We construct a homomorphism from $W_{\infty} / \cap_{k=1}^{\infty} \widetilde{W}^{k}_{}$
to a countable direct product of groups $\mathbb{Z}_2$, i.e., in $\prod\limits_{i=1}^{\infty }{\mathbb{Z}_2}$, which by Theorem 10.13 \cite{Kainth} is homeomorphic to Cantor set.

\begin{theorem}
The order of $\widetilde{A}_{n^j}^{(0,k)}$ is $(\frac{1}{2})^{k+1}\cdot (n!)^{(\frac{n^{(k+1)} - 1}{n - 1})}$.
The order of the generalized alternating group of $m$-th level
$\widetilde{A}_{n^m}^{(m)}$ is $\frac{1}{2}(n!)^{n^{m}}$, for $m:$ $0 \leq m \leq k$.
\end{theorem}
\begin{proof}
By successively applying Lemma \ref{Invunderconj} to the set of permutations from $X^l$, $1<l \leq k$, we see that the condition \ref{evensumk} of Definition \eqref{A^k} holds for each level $X^l$ after conjugation as well as before.
Thus, $\widetilde{A}_{n^j}^{(1,k)}$ is normal in $\mathop{\wr}\limits^{k}_{i=1}  S_{n_i}$.

To find the index of the group $|W_k : \widetilde{A}_{n^j}^{(0,k)}|$, we will count the order of
$W_k$
$$ |S_{n} \wr \ldots \wr S_{n}| = (n!)^{(1+n+n^{2}+\ldots+n^{k})}
= (n!)^{(\frac{n^{k+1} - 1}{n - 1})}.$$

Now we will calculate the order of
$\widetilde{A}_{n}^{(0,k)}$.
For this,
let us count the size of possible permutation tuples satisfying \eqref{evensummk}
that can operate at the $m$ levels $0 \leq m \leq k-1$.
These permutation tuples form subgroups of $W_k(l)$ levels of $Aut{{X}^{[k]}}$. They first appeared in the paper \cite{SkCommEur}, where tuples with an even product of permutations were defined on an arbitrary set of permutations on the vertices of the $k$-th level of $Aut X^{[k]}$.

In order to the number of odd permutations in $X^k$ be even, there must be an even number of $2l$ odd permutations, to satisfy the Definition \eqref{evensumk} and especially the equality \eqref{evensumk}. In order to compute a sum of ${n^{k} \choose 2l}$, consider a sum in the left-hand side of the equality $(1+1)^n=2^n$ and $(1-1)^n=0$, which is $\left({n^{k} \choose 0}+{n^{k} \choose 2}+ {n^{k} \choose 4}+ \ldots \right)$ and divide it by 2

\begin{equation}\label{placement}
\left(\sum_{l=0}^{{n^{k}}}{n^{k} \choose 2l}(1+(-1)^l)\right):2 =
2^{n^{k}-1}.
\end{equation}

Since the symmetric group $S_n$ has the same number of even and odd permutations $\frac{n!}{2}$, then each vertex of the $k$-th level can have one of $\frac{n!}{2}$ permutations, both even and odd. However, the total number of such permutation placements is described above in \eqref{placement}.

Thus, the order of level $m$-th subgroup
$t_{m}:= |G_m| = 2^{n^{m}-1}\times(\frac{n!}{2})^{n^{m}}
= \frac{1}{2}(n!)^{n^{m}}$.

Finally the order of the whole group
$$|\widetilde{A}_{n^j}^{(0,k)}| = \prod_{m=0}^{k-1}|G_m| = (\frac{1}{2})^{k+1}\cdot (n!)^{\frac{n^{(k+1)} - 1}{n - 1}}.$$
\end{proof}

To study the parity of elements at all levels, we factorize by the \textit{normal subgroup} $\tilde{A}_{n^i}^{(0,k)}$.
	\begin{theorem} \label{Isomoflattice}
 The quotient $\mathop{\wr}\limits^{k}_{i=0}S_{n_i}$ by $ \widetilde{A}_{n^j}^{(0,k)}$ is isomorphic to $\prod\limits_{i=0}^{k} \mathbb{Z}_{2}$.
 The order of the quotient  $\left. \mathop{\wr}\limits^{k}_{i=0}  S_{n_i} \right/ \widetilde{A}_{n^j}^{(0,k)}$ is $2^{k+1}$.
\end{theorem}
\begin{proof}
In order to prove, we need to construct homomorphic mappings from $\mathop{\wr}\limits_{i=1}^{{{k}^{}}} S_{n}$ to $\prod\limits_{i=0}^{k} \mathbb{Z}_{2}$ that certify the isomorphism
$$\left. {}^{\mathop{\wr}\limits_{i=0}^{{{k}^{}}}  S_{n}} \right/
{}_{ \widetilde{A}_{n^j}^{(0,k)}} \simeq  \prod_{i=0}^{k} \mathbb{Z}_{2}. $$

Consider the homomorphism $\varphi_l$ from $l$-th level subgroup $W_k(l)$, $l\leq k$ \cite{SkCommEur} onto ${\mathbb{Z}_{2}}$ such that $\varphi_l (\alpha )=\sum\limits_{i=1}^{{{n}^{l}}}{{{dec(\alpha_{li})}}}\bmod 2$ for $\alpha \in W_k$ and as a consequence $\varphi_l(\alpha \cdot \beta )=\varphi_l (\alpha )\circ \varphi_l (\beta )=\sum\limits_{i=1}^{{{n}^{l}}}{{dec({\alpha}_{li})}}\texttt{mod2}+\sum\limits_{i=1}^{{{n}^{l}}}{{dec({\beta}_{li})}}\texttt{mod2}$, where $\alpha \cdot \beta \in W_k$ and ${{{\alpha}_{li}}}$ is a v.p. of an automorphism $\alpha$ state in $v_{li}$. For each $l$ there exists $a_{li} \in S_n \ A_n$, so $\phi(a_{li}) = 1$, which establishes the surjectivity of the mapping.

Note that $\varphi_l(\alpha )$ is an equivariant mapping from $W_k(l)$ to $\mathbb{Z}_{2}$.

To check the main property of a homomorphism we decompose it $$\varphi_l(\alpha \cdot \beta )=\varphi_l (\alpha )\circ \varphi_l (\beta )=(\sum\limits_{i=1}^{{n}^{l}}{{{dec}}(\alpha_{li})}\texttt{mod2}+\sum\limits_{i=1}^{{{n}^{l}}}{{dec({\beta}_{li}}})\texttt{mod2}) \texttt{mod2}.$$
On the left side, we transform a homomorphic image $\varphi_l(\alpha \cdot \beta )$ by the law \eqref{evendec} as
\begin{gather}
\varphi_l(\alpha \cdot \beta ) =  \sum\limits_{i=1}^{{{n}^{l}}}{{{dec}}(((\alpha \cdot \beta)_{li})} \texttt{mod2}\equiv \nonumber \label{conjmonolith}
\\[-3mm]
\label{evendecHom} \nonumber 
   \\[-3mm]
  (\sum\limits_{i=1}^{{{n}^{l}}}{{{dec}}(\alpha_{li})} + \sum\limits_{i=1}^{{{n}^{l}}}{{{dec}}(\beta_{li})}  -
2m) \texttt{mod2} \equiv (\sum\limits_{i=1}^{{{n}^{l}}}{{{dec}}(\alpha_{li})} + \sum\limits_{i=1}^{{{n}^{l}}}{{{dec}}(\beta_{li})})\texttt{mod2}, m \in \mathbb{N}.  \nonumber
\end{gather}

By successively applying Lemma \ref{Invunderconj} to the set of permutations from $X^l$, $1<l \leq k$, we see that the condition \ref{evensumk} from Definition \eqref{A^k} holds for each level $X^l$ after and before a conjugation.
Thus, $\widetilde{A}_{n^j}^{(0,k)}$ is normal in $\mathop{\wr}\limits^{k}_{i=0}  S_{n_i}$. Then the quotient can be computed and it turns out to be isomorphic to $ \mathop{\prod}\limits_{m=1}^{k}\mathbb{Z}_{2}$, furthermore, it is the subgroup of $W_k$.
Therefore, we obtain the order of the quotient.
$$\left.|\mathop{\wr}\limits^{k}_{i=0} S_{n_i} \right/ \widetilde{A}_{n^j}^{(0,k)}| = 2^{k+1} = | \prod_{i=0}^{k}\mathbb{Z}_{2}|.$$
That completes the proof.
\end{proof}

\begin{corollary}
    The normal subgroup lattice of  $\left. { {\mathop{\wr}\limits_{i=0}^{{{k}}}  S_{n}} } \right/
{ \widetilde{A}_{n^k}^{(0,k)}}$
        is isomorphic to the normal subgroup lattice of $\prod\limits_{i=0}^{k} \mathbb{Z}_{2}.$
        The number of normal subgroups in   $\left. \mathop{\wr}\limits^{k}_{i=0}  S_{n_i} \right/ \widetilde{A}_{n^j}^{(0,k)}$ is equal to
$2^{2^{k+1}}.$
\end{corollary}
\begin{proof}
According to Theorem \ref{Isomoflattice} the quotient group $\left. { {\mathop{\wr}\limits_{i=1}^{{{k}}}  S_{n}} } \right/
{ \widetilde{A}_{n^k}^{(0,k)}}$ is isomorphic to $\prod\limits_{i=0}^{k} \mathbb{Z}_{2}.$ We establish correspondence between normal subgroups.
Since there is a bijection between normal subgroups under an isomorphism of groups $\prod\limits_{i=1}^{k} \mathbb{Z}_{2}$, there is also an isomorphism between lattices of normal subgroups.

In view of the fact that each normal subgroup of $\prod\limits_{i=0}^{k} \mathbb{Z}_{2}$ is generated by an arbitrary set of elements of the whole group we deduce that the number of normal subgroups is equal to the Boolean of the set of subgroups of $\prod\limits_{i=1}^{k+1} \mathbb{Z}_{2}$,
$$\sum_{l=0}^{2^{k+1}}{2^{k+1}\choose l} = 2^{2^{k+1}}.$$
According to Theorem \ref{Isomoflattice} lattice of normal subgroups of
   ${}^{W_k}  {/}  {}_{ \widetilde{A}_{n^j}^{(0,k)}}$
is isomorphic to lattice of normal subgroups $\prod\limits_{i=0}^{k} \mathbb{Z}_{2}$, then, the number of normal subgroups in the entire group is equal to $2^{2^{k+1}}$.
\end{proof}

\begin{remark}\label{n^3}
The order of $\widetilde{A}_{n^2}^{(2)}$ is $(n!)^{n^2}\cdot 2^{-1}$.
\end{remark}
\begin{proof}
The order of a $n^2$-fold direct product of permutation groups $S_n$ is $(n!)^{n^2}$.
Exactly half of the permutation tuples from the $n^2$-th level have an even sum of decrements, which determines the order of the group.
\end{proof}

\begin{theorem}
    Let $H < W_k$, if $A_{{{n}^{m}}}^{(k)}<H$, where $d(H)=l$ and $l:\,\,k-l=m$ and a complement kernel $C$ for $A_{{{n}^{m}}}^{(k)}$ to $H$ is invariant in $W_k$, then $C \wr \hspace{-1.5 mm}\wr A_{{{n}^{m}}}^{(k)} \lhd W_k $.
    \end{theorem}
\begin{proof}
Since a subwreath product $C \wr \hspace{-1.1 mm}\wr A_{{{n}^{m}}}^{(k)} \lhd W_k$ is a subgroup of iterated wreath product and thus a split extension $C \ltimes A_{{{n}^{m}}}^{(k)}$, then the conjugation of elements from each factor occurs level by level with respect of top subgroup action, so in conditions of our theorem it remains to verify that such conjugation preserves a subdirect product $A_{{{n}^{m}}}^{(k)}$ on last level.
Since in this case $A_{{{n}^{m}}}^{(k)}$ is base, all automorphisms of form $g_{(v_{li})}|_{X^{[m]}} \in S{{t}_{W}}\left( l  \right) $ permute v.p. inside of a correspondent tuple of vertexes of $A_{{{n}^{m}}}^{(k)} $ in the base subgroup  of  $Aut{{v}_{li}}{{X}^{[m]}} $  containing  $S{{t}_{W}}\left( {{v}_{li}} \right)$. This means that v.p. from $A_{{{n}^{m}}}^{(k)} $ are rearranged only with each other by automorphisms  from  $W_{{v}_{li}}$, therefore $\sum\limits_{i=k-l+1}^{m}{dec\left( {{a}_{ki}} \right)=2t}$.
Thus, the restriction $g_{(v_{li})}|_{X^{[k-l]}}$ of the action of an automorphism $g\in AutX^k$ on the subtree $v_{li}X^{[k-l]}$ keeps the level subgroup $A_{{{n}^{m}}}^{(k)}$ invariant.

Taking into account that each element of $H$ acts by conjugation non-primitively on blocks (that is, it permutes the sets of vertices of blocks without breaking these blocks) of $n^m$ elements of $X^k$ forming blocks of imprimitivity of size $n^m$ we see that condition \eqref{evensummk} of Definition \ref{A^mk} on blocks with limits $[i=(s-1)n^m+ 1, sn^m]$ is satisfied, because each block is the support of the level subgroup of the vertex stabilizer $W_{{v}_{li}}$
for the correspondent $i$. The number of such blocks is $n^{k-m}$. The top subgroup $C$ is normal by the condition.

$A_{{{n}^{m}}}^{(k)}$ acts trivially on vertices of levels $X^j$ $j<l$ conjugation of elements in form of wreath recursion with root permutation on $l$ level takes form
\begin{gather}
{{\pi }^{-1}}\left( a_{{{\pi }^{-1}}( (s-1)n^m+ 1)}^{-1},...,a_{{{\pi }^{-1}}( n^m)}^{-1} \right) \, \sigma \, ({{b}_{(s-1)n^m+ 1}},...,{b}_{( n^{m})})\pi \left({{a}_{(s-1)n^m+ 1}},...,{{a}_{(n^{m})}} \right)= \nonumber
\\[-3mm]
\label{conjXln{k-l}}
   \\[-3mm]
\pi^{-1}  \sigma {{\pi }} \left( a_{{{\pi }^{-1}}\sigma \pi (( (s-1)n^m+ 1))}^{-1}{{b}_{\pi (1)}}{{a}_{((s-1)n^m+ 1)}},...,a_{{{\pi }^{-1}}\sigma \pi (( sn^m))}^{-1}{{b}_{\pi (sn^m)}}{{a}_{(sn^m)}} \right).  \nonumber
\end{gather}
for each level subgroup $v_{l,i}X^{[m]}$.
Moreover, in this case we can reformulate the statement if $  S{{t}_{W}}\left( {{v}_{li}} \right)  \triangleleft H$, $d(H)=l$ and  $l:\,\,k-l=m$, then  $\widetilde{A}_{{{n}^{m}}}^{(k)}\triangleleft S{{t}_{W}}\left( {{v}_{li}} \right )$. In other words, in this case we have transitivity of normality.
\end{proof}
\begin{example}
    For instance $E \wr E \wr \hspace{-1.5mm}\wr \widetilde{A}_{{{n}^{2}}}^{(2)} \wr \hspace{-1.5mm}\wr \widetilde{A}_{{{n}^{}}}^{(3)}$ is normal in $S_n \wr S_n \wr S_n \wr S_n$.
\end{example}

\subsection{Examples and computations in $S_{n} \wr S_{m}$}
As is well known $S_2 \wr S_2  \simeq D_4$ and has 6 normal subgroups, of which 4 are proper.
Computations in GAP lead us to following results and conclusions:

    The group $S_{3} \wr S_{3}$ has 10 normal subgroups,  where  $G_0=E$  and $G_1=S_{3}\wr S_{3}=\langle (1,2,3), (1,2),
    (4,5,6), (4,5), (7,8,9), (7,8), (1,4,7)(2,5,8)(3,6,9), (1,4)(2,5)(3,6)\rangle$.

In total we have 8 proper normal subgroups in $S_{3} \wr S_{3}$.

By analyzing their structures in terms of a semidirect product and determination them as the subgroup of wreath product with
respect of their orders we obtain:

    1) $S_3 \wr \widetilde{A_3}  \simeq  (((C_3 \times C_3 \times C_3) \rtimes (C_2 \times C_2)) \rtimes C_3) \rtimes C_2,  $
    here are 3 non equal isomorphic subgroups, also we compute $ord(S_3 \wr  \tilde{A_3}  )=((6^3):2) \cdot 6=3^4 \cdot
    2^3=108\cdot 6=648$.

    2) $\langle  (1,2,3), (1,2), (4,5,6), (4,5), (7,8,9), (7,8), (1,4,7)(2,5,8)(3,6,9), (1,4)(2,5)(3,6)  \rangle=S_3\wr S_3$.

All subgroups, except $G_7$, from this list can be easily identified and presented with new indexing.

 We denote $S_3\wr S_3$ by $H_1$ and we compute that $\left[ H_1 =  (((( C_3 \times  C_3 \times C_3) \rtimes (C_2 \times C_2)) \rtimes
 C_3) \right. \\ \left. \rtimes C_2)  \rtimes C_2, 1 \right].$ Order of $H_1$ is 1296. \\
We find the structures of these 8 proper normal subgroups of $S_3\wr S_3$:\\
$H_2= [(((C_3 \times C_3 \times C_3) \rtimes (C_2 \times C_2)) \rtimes C_3) \rtimes C_2, 3],$ \\
 $H_3= [((C_3 \times C_3 \times C_3) \rtimes (C_2 \times C_2)) \rtimes C_3, 1],$ \\
$H_4=[ (C_3 \times C_3 \times C_3) \rtimes (C_2 \times C_2), 1 ],$\\
$H_5=[ (C_3 \times C_3 \times C_3) \rtimes C_2, 1 ],$ \\
$H_6=[ C_3 \times C_3 \times C_3, 1 ],$ \\
$H_7=[ S_3 \times S_3 \times S_3, 1 ] ].$

\subsection{Structural theorems}
Thus, $S_{3} \wr S_{3}$ contains 8 proper normal subgroups.
Three of them  having the same orders 648 (denoted by $H_2$) are determined by semidirect products of the same groups but with different homomorphisms in the $Aut(Ker(H_2))$,
therefore they not permutation isomorphic).
 Furthermore a form of their generators we can deduce the splitability
of this groups.

    1) Consider the three isomorphic subgroups corresponding to $G_7$ having GAP structure description  $H_2 = ["(((C_3 \times C_3 \times C_3) : (C_2 \times C_2)) : C_3) : C_2", 3 ]$ here are 3 non equal non isomorphic subgroups of order 648. Algebraic structure of them is the following:
\[H_2  \simeq \left( \left( \left( C_3\text{ }\times \text{ }C_3\text{ }\times \text{ }C_3 \right)\rtimes \left( C_2\text{ }\times
\text{ }C_2 \right) \right)\rtimes \text{ }C_3 \right)\text{ }\rtimes \text{ }C_2,\]
 these subgroups have order 648 and identifiers according to GAP system of this small group: $IdGroup =703, 704$ and $705$ respectively.

  We will denote these 3 subgroups by ${{H}_{703}},\,\,\,{{H}_{704}}$ and ${{H}_{705}}$ correspondingly to their "IdGroup" in GAP.
    At the same time, we have 2 different groups isomorphic up to semidirect product to $H_2$ with this order $ ord(S_3 \wr \widetilde{A_3} )= ord(A_3 \wr S_3) =((6^3):2) \cdot 6=648$.

   The commutator subgroups of $H{{}_{703}}$ and $H{{}_{704}}$ have order 324 and the same GAP
  identifier $"IdGroup(H'{{}_{703}}) = IdGroup(H'{{}_{704}}) = 160"$. Thus $H{{}_{703}}$ and $H{{}_{704}}$ have the same multipliers in the semidirect product with non equal actions, and they not equal because of different embeddings in $S_3 \wr S_3$. Note that $ ord(S_3 \wr \widetilde{A_3} )' =324$ because of $(S_3 \wr \widetilde{A_3} )' \simeq  A_3 \wr \widetilde{A_3} $.
   This means that commutator of these subgroups is the same therefore
  ${{H}_{703}}$ and ${{H}_{704}}$ \textbf{are isomorphic but they are embedded} in ${{S}_{3}}\wr {{S}_{3}}$ in \textbf{different
  ways} as different copies. The center of these subgroups are trivial subgroup. For the sake of clarity, we will \textbf{denote the second copy} of $H_{704}$ as $\widetilde{{\mathbb{B}}}_n$.
  The commutator subgroup of
  \[{{S}_{3}}\wr {{\widetilde{A}}_{3}} \simeq H_{703} \simeq  H_{704} \]
   is exactly ${{A}_{3}}\wr {{\widetilde{A}}_{3}}$, so the top subgroup has order in 2 times less then order of ${{S}_{3}}$.

 Furthermore, the normal subgroup $A_{3} \wr S_{3}$ has the commutator subgroup $e \wr \widetilde{A}_{3}$ with the same
 order 108, which is in 3 times lesser than $ ord(S_3 \wr \widetilde{A_3} )' =324$. Thus, the \textbf{third copy} of $H_2$ from the GAP list is  $A_{3} \wr S_{3}=H_{705}.$

The invariant subgroup lattice of $S_n \wr S_n$ for cases $n\equiv 1\, (\mathrm{mod 2})$ and $n\equiv 0\, (\mathrm{mod 2})$ are presented on Fig. 1. and Fig. 2.
 \begin{figure}[h]
    \centering
    \begin{minipage}{0.48\textwidth}
        \centering
        \includegraphics[width=\linewidth]{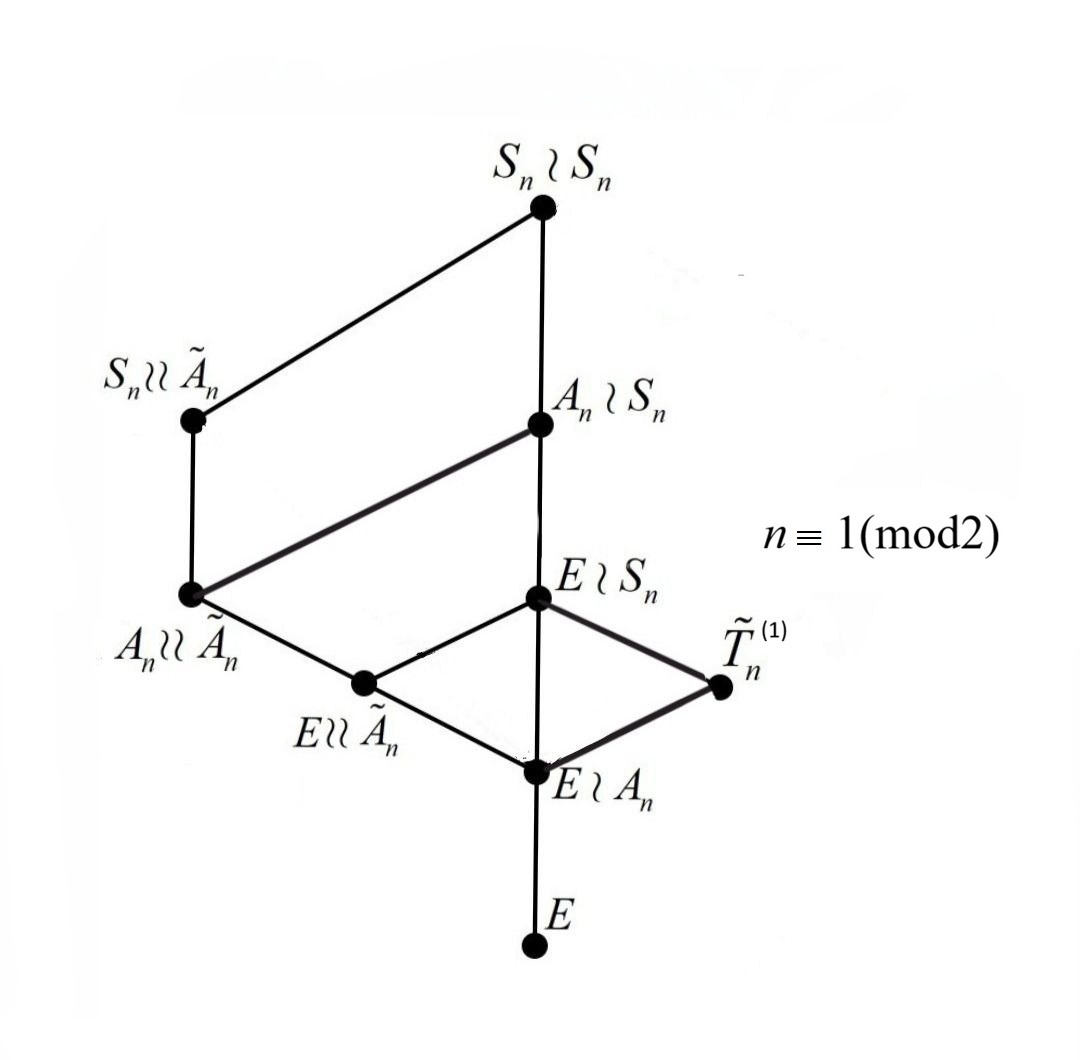}
        \caption{Lattice of invariant subgroups $S_n \wr S_n$ for the case $n\equiv 1\, (\mathrm{mod 2})$}
    \end{minipage}
    \hfill
    \begin{minipage}{0.48\textwidth}
        \centering
        \includegraphics[width=\linewidth]{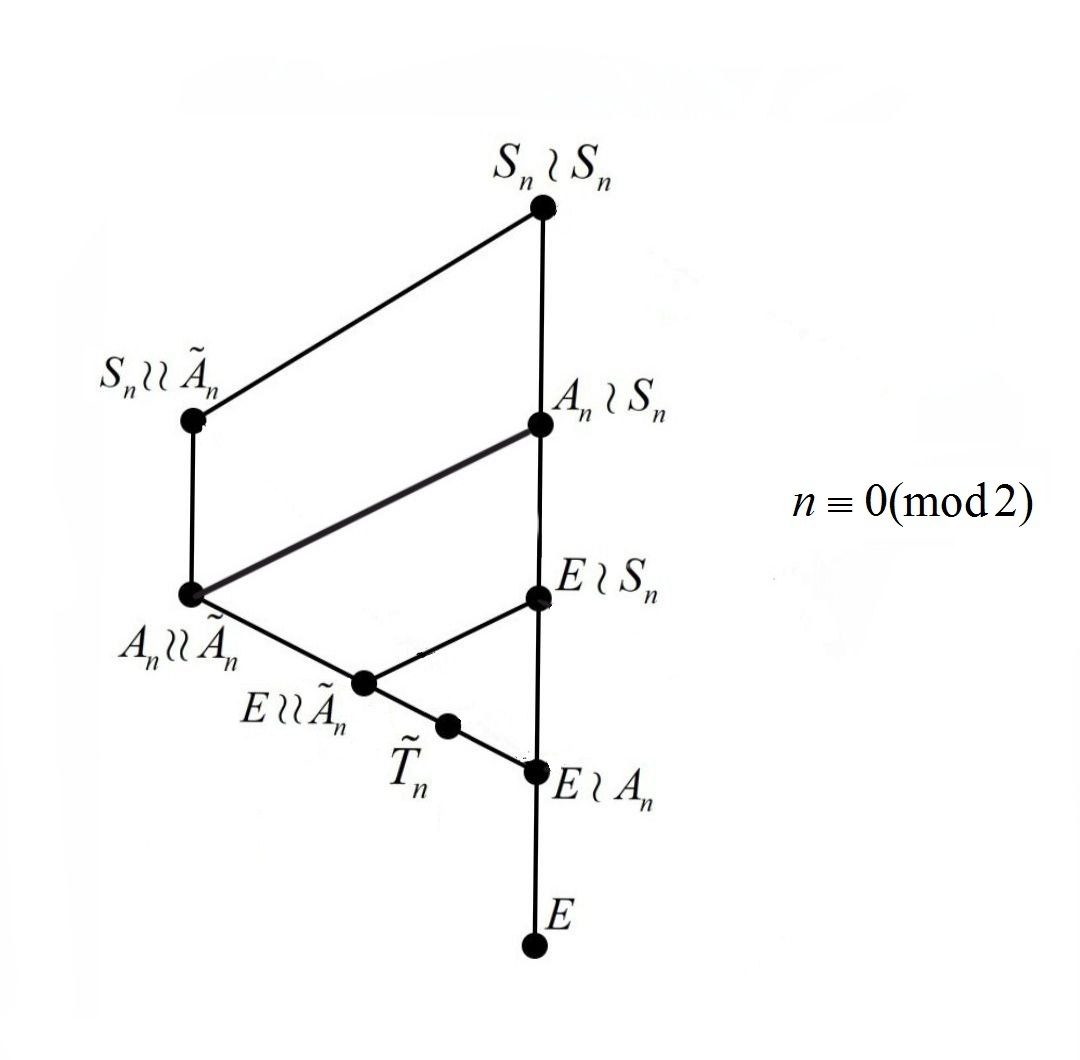}
        \caption{Lattice of invariant subgroups $S_n \wr S_n$ for the case $n\equiv 0\, (\mathrm{mod 2})$}
    \end{minipage}
\end{figure}
According to the statement of Corollary\ref{Mon Short E wr A_n} subgroup $E\wr A_n$ is the atom of each of these lattices.

   Based on the subgroups computed above and an analysis of their properties, we formulate a theorem.
\begin{theorem} \label{Thn=3}
Proper normal subgroups in $S_{3} \wr S_{3}$ are of the following types:

\begin{enumerate}
\item subgroups that act only on the second level (stabilizing the first level \cite{SkCommEur, GrNe}) are
$$ E \wr\hspace{-1mm}\wr \widetilde{A}_{3}, \, \widetilde{T_{3}}, \, E \wr
S_{3},  E \wr A_3,$$

\item subgroups that act on both levels are
$A_3 \wr\hspace{-1mm}\wr \widetilde{A}_{3}, \,    \, S_3 \wr\hspace{-1mm}\wr \widetilde{A}_{3},  S_3 \wr\hspace{-1mm}\wr \widetilde{B}_{3}, \, A_3 \wr S_{3}, $
\end{enumerate}
wherein the subgroup $S_3 \wr \widetilde{A}_{3} \simeq S_3 \rightthreetimes (\underbrace{ S_3 \boxtimes S_3 \boxtimes S_3 }_3 )$ endowed with the subdirect product satisfying condition \eqref{evensum}. 

In total, there are 10 normal subgroups, including 8 of their proper normal subgroups in $S_{3} \wr S_{3}$.
\end{theorem}
\begin{proof}

    The group $W_3$ has exactly eight normal subgroups, by virtue of the classification theorems proved above and the analysis of the number of isomorphic copies of normal subgroup in the GAP system presented in this section.
\end{proof}

\begin{theorem} \label{Thn}
Proper normal subgroups in $S_{n} \wr S_{m}$, where $n, m
\ge 3$ with $n, m \neq 4$ are of the following types:

\begin{enumerate}
\item subgroups that act only on the second level (stabilizing the first level \cite{SkCommEur, GrNe}) are
$$ E \wr\hspace{-1mm}\wr \widetilde{A}_{m}, \, \widetilde{T_{m}}, \, E \wr
S_{m},  E \wr A_m,$$

\item subgroups that act on both levels are
$A_n \wr\hspace{-1mm}\wr \widetilde{A}_{m}, \,    \, S_n \wr\hspace{-1mm}\wr \widetilde{A}_{m},  S_n \wr\hspace{-1mm}\wr \widetilde{B}_{m}, \, A_n \wr S_{m}, $
\end{enumerate}
wherein the subgroup $S_n \wr \widetilde{A}_{m} \simeq S_n \rightthreetimes (\underbrace{ S_m \boxtimes S_m \boxtimes S_m \boxtimes \ldots \boxtimes S_m}_n )$ endowed with the subdirect product satisfying condition \eqref{evensum}.

In total, there are 10 normal subgroups, including 8 of their proper normal subgroups in $S_{n} \wr S_{m}$.
\end{theorem}

\begin{proof}

Firstly, we remark that $E \wr S_{m}$ is normal as the base \cite{ShNorm, SkCommEur}  of $S_n \wr S_m$.

In view of Proposition \ref{Mon Short} $E \wr A_m$ is minimal non-trivial normal subgroup. This implies that $ S_n \wr S_ {m} $ is \emph{monolithic} group and its \emph{monolith} is $ e \wr A_ {m} $.

The normality of $ e \wr \widetilde{A}_ {m} $ is justified in Proposition \ref{A2subgr}.

In order to prove the normality of $S_{n} \wr \widetilde{A}_{m}$
it is sufficient to validate the coincidences of the normal closure in $S_{n} \wr S_{m}$ of its generators to the correspondent isomorphic copy of the initial
subgroup. As generators of $S_{n} \wr \widetilde{A}_{m}$ may be or $\widetilde{a} = [e]_{1}, \ [(i_1, i_2, i_3), e,  \ldots ,e]_{2}$ either $\widetilde{b} = [e]_{1}, \ [(i_1, i_2), (i_4, i_3), e,  \ldots ,e]_{2}$, $1 \leq i_j \leq n$ and necessary generator of $S_n \wr E$ is $\widetilde{s_j} = [(j_1, j_2)]_{1}, \ [e, e, e,  \ldots ,e]_{2}$, $(j_1 \neq j_2)$. A normal closure of generator $(j_1, j_2)$ in $S_n$ is $S_n$.
Recall the normal closure of \emph{generators of} $e \wr \widetilde{A}_{m}$ is found in Proposition \ref{A2subgr} and it coincide with $e \wr \widetilde{A}_{m}$.

Observe that the product of powers $\widetilde{a}$ and $\widetilde{b}$ satisfies the condition \ref{evensum}, thence such products belong to $ E \wr \widetilde{A}_{m}$.

Constructing a normal closure of $\left \langle \widetilde{a}, \widetilde{b}, \widetilde{s_j} \right \rangle$ and taking into account Lemma \ref{Invunderconj} about invariant of even decrement sum under conjugation provided $j=1$ we see the sum of decrements in the right part of \eqref{conjXk} complies with \eqref{evensum}.
Thus, $S_{n} \wr \widetilde{A}_{m}$ is normal.

As mentioned above, $A_{n}$ is generated by all $3$-cycles and the fact that all 3-cycles can be generated by conjugating of
$(i_{1}i_{2}i_{3})$ by even permutations. Therefore normal closure of additional generator $[(i_{1}i_{2}i_{3})]_{1}, [e, \ldots,
e]_{2}, \
i_{1}, i_{2}, i_{3} \in \{1,2,\ldots,n\}$ is  $A_{n} \wr E$.

Normality of $A_{n} \wr S_{m}$ is obvious in view of $A_n \lhd S_n$ and isomorphism of $E \wr S_{m}$ to the base of $S_{n} \wr S_{m}$.

2.  Here we prove that the groups mentioned in item 2 are normal with using the fact of normality in $S_{m}$ of subgroups
appearing on second level. Such groups are normal in $S_{m}$.

a)
Normality of $A_{n} \wr \widetilde{A}_{m}$ immediately follows from  $A_{n} \wr \widetilde{A}_{m} = (S_{n} \wr \widetilde{S}_{m})'$,
because of the defining condition \eqref{evensum} of $\widetilde{A}^{(1)}_{m}$ accords with condition of commutator subgroup \cite{Meld, SkCommEur}, satisfying the condition \eqref{evensum} is equivalent to the parity of the product $\prod \limits_{i=1}^n\pi_{1i} \in A_m$ of permutations from the first level, which characterizes the commutator subgroup.
According to \cite{SkCommEur, Gural_2010, Meld} this subgroup is a commutator subgroup of $S_{n} \wr S_{m}$ hence it is normal.

For additional goals we find {generators whose normal closure is $A_{n} \wr \widetilde{A}_{m}$}. These are the same as for
 $e \wr \widetilde{A}^{(1)}_{m}$ and one additional generator $\widetilde{s} = [(i_{1},i_{2},i_{3})]_{0}, [e, \ldots, e]_{1}, \ 1\leq i_{1}, i_{2}, i_{3} \leq n$ completing it to $A_{n} \wr \widetilde{A}^{(1)}_{m}$.
  Therefore the normal closure of $\widetilde{s}, \widetilde{a}, \widetilde{b}$ is
$A_{n} \wr \widetilde{A}_{m}$.

b) \emph{Generating elements of} normal closure equal to $S_{n} \wr
\widetilde{A}_{m}$ have the same elements
as for $A_{n} \wr \widetilde{A}_{m}$ except one additional generator presented by the next tableau:
$$t= [(i_{1},i_{2})]_{1}, [e, \ldots, e]_{2}, \
i_{1}, i_{2} \in \{1,2,\ldots,n\}, (i_1 \neq i_2).$$

The last generator completes normal closure of active group to $S_{n} \wr
E$. Therefore the normal closure of $t,a,b$ is
$S_{n} \wr
\widetilde{A}_{m}$.

At least the normality of $T^{(1)}_{n}$ is established in Proposition \ref{T^{(1)}_{n}}.
\end{proof}

\subsection{Future researches}
The normal rank \cite{Dash} of $W_k$ invariant subgroups will be studied in our next paper.

\textbf{Conclusion.}The monolith of these wreath products has been
investigated.
We have shown that $S_{n_1}\wr S_{n_2} \wr S_{n_3}$, $n, m \in
\mathbb{N}$ and $ S_n \wr S_ {m} $ are the monolithic groups and its \emph{monoliths} are $ e \wr e \wr A_ {m}$ $ e \wr A_ {m}$. We have found normal subgroups and their properties for finite and infinite iterated wreath products $S_{n_1}\wr S_{n_2} \wr S_{n_3}$, $n, m \in
\mathbb{N}$. All classes of normal subgroups in infinite wreath product of symmetric groups have been studied up to isomorphism. The topology based on normal subgroups in the infinite iterated wreath product of permutation groups has been investigated.

The author declares no conflict of interest.

{\scriptsize \textbf{Sources of Funding for Research Presented in a Scientific Article. 
This work was partially supported by a grant from the Simons Foundation (SFI-PD-Ukraine-00017674, Ruslan Skuratovskii).
}}


\begin{thebibliography}{9}

\bibitem {SkuNormNPU} \emph{ Skuratovskii R.V.,} Invariant structures of wreath product of symmetric groups. Naukovuy Chasopus of Science hour writing of the NPU named after M.P. Dragomanova. (in ukrainian) Series 1. Physics and Mathematics. 2009. Issue 10. P. 163-178.



\bibitem{ShLav}  \emph{Lavrenyuk, Y.V., Sushchanskii, V.I.} Lattice of invariant subgroups of a group of local isometries of the boundary of a spherically homogeneous tree (in russian). Ukr Math J  60, 1574-1580 (2008). https://doi.org/10.1007/s11253-009-0154-8

\bibitem{LeshCong} Yurij Leshchenko. (Infinitely Iterated Wreath Product of Elementary Abelian Groups) Ukrainian Mathematical Congress - 2009. Source: [https://www.imath.kiev.ua/~congress2009/Abstracts/Leshchenko.pdf]

\bibitem{LeshPap} Yurij Leshchenko. Infinitely iterated wreath power of elementary Abelian groups. Matematychni Studii. V. 31, No.1.

\bibitem{ShNorm}
\emph{ Sushchansky  V.~I.}
\newblock Normal structure of the isometric group of metric spaces of $p$-adic integers.
\newblock {\em Algebraic structures and their application. Kiev}, Visn. of KNU, 1988. pp.
  113-121.

 \bibitem{ShLav93} V.I. Sushchanskii, Normal Structure of Isometry Groups of Semifinite Baire Metrics. Infinite Groups and Related Algebraic Structures
[in Russian], Institute of Mathematics, Ukrainian National Academy of Sciences, Kiev (1993).

\bibitem {aviv}  \emph{Aviv Rotbart,} Generator sets for the alternating group. Seminaire Lotharingien de Combinatoire 65 (2011), Article.

\bibitem {SubAlg} \emph{ Birkhoff, Garrett}  (1944), "Subdirect unions in universal algebra", Bulletin of the American Mathematical Society, 50 (10): 764-768, doi:10.1090/S0002-9904-1944-08235-9, ISSN 0002-9904, MR 0010542.


\bibitem{Meld}  \emph{J.D.P. Meldrum,} Wreath Products of Groups and Semigroups.  Pitman Monographs and Surveys in Pure and
    Applied Mathematic. 1st Edition. Jun (1995). 425 p.

\bibitem{LB} \emph{Leonov U. G.}, On the representation of groups approximated by finite $p$-groups.
//{\em  Ukr. Math. Journ.}, 2011. -- {\em V.} 63, {\em No.} 11.--{\em P}.
1518-1511.


\bibitem{Ne}  \emph{Nekrashevych V.,} Self-similar groups. International University Bremen. American Mathematical Society.
    Monographs. Volume 117. 230 p.

\bibitem{Lav}  \emph{Lavrenyuk Y.,} On the finite state automorphism group of a rooted tree. Algebra and Discrete Mathematics
Number 1. (2002). pp. 79-87.






\bibitem{Kal}   \emph{ Kaloujnine~L.\,A.} Sur les $p$-group de Sylow. 
    // C.~R.~Acad. Sci. Paris. --- 1945. --- {\bf 221}. --- P.\,222--224.

\bibitem{Bez}   \emph{Bezuschak O. O.}, Splittable normal subgroups of isometry group of a generalized metric space, Mathematychni Studii, 17 (2002) pp. 29-40.


\bibitem{SkRendi} \emph{ Skuratovskii R. V.,  Williams A.} "Irreducible bases and subgroups of a wreath product in applying to diffeomorphism groups acting on the Mobius band", \textbf{2021}.
    Rendiconti del Circolo Matematico di Palermo Series 2, 70(2), 721-739. https://doi.org/10.1007/s12215-020-00514-5.

\bibitem{Gural_2010} \emph{ Guralnick R.,} Commutators and wreath products  // Contemporary Mathematics. Volume 524, 2010.


\bibitem{Masl} \emph{ N. V. Maslova, D. O. Revin}, On the Pronormality of Subgroups of Odd Index in Some Direct Products of
Finite Groups, \emph{Journal of Algebra and Its Applications}, 22:04 (2023), 2350083, doi: https://doi.org/10.1142/S0219498823500834

\bibitem{DrSku}
\emph{Drozd, Y.A., Skuratovskii R.V.,} Generators and relations for wreath products of groups. {\it Ukr. Math. J.}  (2008), {\em
60}, pp. 1168--1171.


\bibitem{Sachk}
\emph{ Sachkov, V.N.},  Combinatorial methods in discrete Mathematics. Encyclopedia of mathematics and its applications 55. Cambridge Press. 2008. P. 305.

\bibitem{Dash}
\emph{Dashkova O. Yu.} On groups of finite normal rank. \small{Algebra Discrete
Math.} 2002. 1,  No. 1. P. 64-68.

\bibitem{Bidw}
\emph{Bidwell, J.N.S.} Automorphisms of direct products of finite groups II.
\small{ Arch. Math.} 91, pp. 111-121 (2008).
https://doi.org/10.1007/s00013-008-2653-5


\bibitem{SkKPI} \emph{Skuratovskii~R.,} Generators and relations for sylows $p$-subgroup of group $S_n$. \emph{Naukovi Visti
    KPI.} \textbf{4} (2013), pp. 94--105. (in Ukrainian)





\bibitem{Dm}
\emph {Dmitruk~U., Suschansky~V.,}  Structure of 2-sylow subgroup of alternating group and normalizers of symmetric and
alternating group. UMJ. (1981), N. 3,  pp. 304-312.

\bibitem{DWard} \emph { Ward D.,} Topics in Finite Groups: Homology Groups,
Pi-product Graphs, Wreath Products and
Cuspidal Characters. Manchester Institute for Mathematical Sciences
School of Mathematics. July (2015) P. 253.

\bibitem{Bir}
\emph{Birkhoff G.} Theory of lattice. M. Nauka 1984, P. 564.

\bibitem{Cox}
\emph{Coxeter G.S., Moser U.O.} Generating elements and determining relations of descrete groups. M. Nauka., 1980, 240p.


\bibitem{SyOp}
\emph{Sushchansky V.I., Sikora V.S.} Operations on the groups of permutations. Chernivtsi: Ruta, 2003.14.



\bibitem{SkCommEur}
\emph{ Ruslan V. Skuratovskii.}
On commutator subgroups of Sylow 2-subgroups of the alternating group, and the commutator width in wreath products. / Ruslan V. Skuratovskii //
European Journal of Mathematics. -- 2021. --  vol. 7, no. 1. -- P. 353-373. https://doi.org/10.1007/s40879-020-00418-9

\bibitem{GrNe}
\emph {R.~ Grigorchuk, V.~Nekrashevich, V.~Sushchanskii,} Automata, Dynamical Systems, and Groups, Trudy mat. inst. imeny Steklova.  (2000), Vol. 231, P. 134--214.







\bibitem{Kolm}
 A. N. Kolmogorov, S. V. Fomin. Elements of the Theory of Functions and Functional Analysis (Dover Books on Mathematics) Dover Books on Mathematics Edition. 1999, 288 pages.

\bibitem{GrBranch}
R. I. Grigorchuk, On branch and just infinite groups, Abstracts of International Algebra Conference in Memory of Kurosh, Moscow, 1998.

\bibitem{GrInfBranch} Grigorchuk, R.I. (2000). Just Infinite Branch Groups. In: du Sautoy, M., Segal, D., Shalev, A. (eds) New Horizons in pro-p Groups. Progress in Mathematics, v. 184. Birkhauser, Boston, MA. $https://doi.org/10.1007/978-1-4612-1380-2-4$

\bibitem{DueselldorfBranch} Alejandra Garrido. On the Congruence Subgroup Problem for Branch Groups   University of Oxford
https://www.math.uni-duesseldorf.de/~garrido/neuchatel14.pdf

\bibitem{Inverselimit}
Jan-Erik Roos. Derived functors of inverse limits revisited // J. London Math. Soc..  V. 73, (2006). no. 1.  pp. 65-83.  doi:10.1112/S0024610705022416

\bibitem{Shaw} Matthew Shaw. Characterization of Cantor Spaces
Matthew Shaw: "Characterization of Cantor Spaces" November 2019, $https://web.math.utk.edu/~freire/teaching/m467f19/Cantor_Spaces_Topology.pdf$


\bibitem{AlgCantor}  Reznichenko E. Algebraic structures on the Cantor set. Jan. 2022,
 (2022). Algebraic structures on the Cantor set. ArXiv. /abs/2207.01003

\bibitem{Kainth}  \emph{Surinder Pal Singh Kainth.} A Comprehensive Textbook on Metric Spaces  Publisher
Springer Singapore (2024), P. 344, DOI https://doi.org/10.1007/978-981-99-2738-8

\bibitem{tails}  Sidney A. Morris. Topology without tears.  Version of February 20, 2012. 1 c Copyright 1985-2012.

\bibitem{Aschbacher} \emph{ M. Aschbacher.}
Finite Group Theory. August 2, 2010 Cambridge Studies in Advanced Mathematics 10)   2nd Edition, 318 pages.
\end{thebibliography}
\end{document}